# Low Mach number limit for the diffusion approximation model in radiation hydrodynamics at equilibrium-diffusion regime

Kwang-Il Choe[1], Dae-Won Choe [1], Myong Chol Pak [2,*]

[1] School of Mathematics, University of Mechanical Engineering Pyongyang, Pyongyang, Democratic People's Republic of Korea

[2] Department of Physics, **Kim Il Sung** University, Taesong District, Pyongyang, Democratic People's Republic of Korea

**Abstract**

The low Mach number limit for the compressible viscous diffusion approximation model arising in radiation hydrodynamics is rigorously justified. For the 3-D Cauchy problem, the solutions in an equilibrium diffusion regime are shown to converge to the solutions of an incompressible Navier-Stokes equations locally and globally in time as Mach number goes to zero, when the effect of the small temperature variation upon the limit is taken into account.

**MSC 2020:** 76N05, 76N10, 35B40, 35A01.

**Keywords:** radiation hydrodynamics, diffusion approximation model, low Mach limit, equilibrium-diffusion regime

## 1. Introduction

The diffusion approximation (also are called the Eddington approximation) model in radiation hydrodynamics can be written in terms of Euler coordinates as follows (see Appendix in [7]):

$$\begin{aligned} &\rho_t + \operatorname{div}(\rho\mathbf{u}) = 0, \\ &\rho\mathbf{u}_t + \rho(\mathbf{u}\cdot\nabla)\mathbf{u} + \nabla P = \mu\Delta\mathbf{u} + (\mu+\lambda)\nabla\operatorname{div}\mathbf{u}, \\ &\rho e_t + \rho\mathbf{u}\cdot\nabla e + P\operatorname{div}\mathbf{u} = \kappa\Delta\theta + 2\mu\mathbf{D}(\mathbf{u}):\mathbf{D}(\mathbf{u}) + \lambda(\operatorname{div}\mathbf{u})^2 - \tilde{\sigma}\theta^4 + \sigma_a n, \\ &\frac{1}{c}n_t - \nu\Delta n = \tilde{\sigma}\theta^4 - \sigma_a n \end{aligned} \tag{1.1}$$

for $\mathbf{x} = (x_1, x_2, x_3) \in \mathbb{R}^3$ and $t > 0$ .Here, the unknown functions $(\rho, \mathbf{u}, \theta, n)$ represent the density, the velocity, the temperature and the radiation field, respectively the pressure $P=P(\rho,\theta)$ and inertial energy $e=e(\rho,\theta)$ are the smooth functions of $\rho$ and $\theta$: $\mu$ and $\lambda$ are the constant viscosity coefficients, $\mu>0, 3\lambda+2\mu\geq 0$: $\kappa>0$ is the heat-conducting coefficient: $c>0$ is the light speed: $\nu>0$ is the diffusion coefficient: $\sigma_a>0$ is absorption coefficient: $\tilde{\sigma}$ is the positive constant defined by

$$0<\tilde{\sigma}=4\pi\sigma_a\alpha\beta^{-4}\int_0^\infty \frac{s^3}{e^s-1}ds<\infty$$

with positive constants $\alpha$ and $\beta$, and $\mathbf{D}=\mathbf{D}(\mathbf{u})$ is the deformation tensor

$$D_{ij}=\frac{1}{2}\left(\frac{\partial u_j}{\partial x_i}+\frac{\partial u_i}{\partial x_j}\right) \quad \text{and} \quad \mathbf{D}(\mathbf{u}):\mathbf{D}(\mathbf{u})=\sum_{i,j=1}^{3}D_{ij}^2.$$

*Corresponding author: mc.pak1001@ryongnamsan.edu.kp

The system (1.1), which consists of the conservation laws of mass, momentum, energy, and the radiative transfer equation, describe a non-equilibrium regime where the state of the radiation is determined by the transport equation $(1.1)_4$. Since its physical importance, complexity, rich phenomena and mathematical challenges, there are huge literatures on the studies of radiation hydrodynamics from the mathematical/physical point of view. For the diffusion approximation model (1.1), the global existence of solutions was proved in Jiang-Xie-Zhang [13,6] and Jiang [7] in the case of the 1-D initial-boundary value problem and 3-D Cauchy problem, respectively. For a non-equilibrium diffusion approximation model in thermally radiative magnetohydrodynamics, which is a generalization of the system (1.1), the global existence of solutions was proved in Jiang [8] and Jiang-Yu [9] for 1-D and 3-D case, respectively. For 3-D Cauchy problem of the system (1.1) with $\mu=\lambda=\kappa=0$, the local existence of solutions was proved in Jiang-Zhou [10]. Recently, Kim-Hong-Kim [14] showed the global existence and decay rates of the smooth solutions for 3-D Cauchy problem of the system (1.1) in the case of the more general gases including ideal polytropic gas.

On the other hand, the diffusion limits in radiation hydrodynamics consist of supposing that one of the transport coefficients is small while the other is large. In particular, when the scattering effect is assumed to be small, the corresponding regime is called an equilibrium diffusion regime. On the contrary, when the scattering effect is dominant, the corresponding regime is called a nonequilibrium diffusion regime.

The diffusion approximation model (1.1), which is valid for optically thick regions where the photons emitted by the gas have a high probability of reabsorption within the region, is particularly accurate if the specific intensity of radiation is almost isotropic ([16]).

The present paper aims to justify the low Mach number limit at equilibrium-diffusion regime for the diffusion approximation model (1.1) when the effect of the small temperature variation upon the limit is taken into account. Indeed, there are some more works on the singular limits for NSFR(the so-called Navier-Stokes-Fourier-Radiation) model. Ducomet-Necasova [1,2] investigated the low Mach number limit for the NSFR system, and Jiang-Li-Xie [12] and Fan-Li-Nakamura [3] studied the nonrelativistic and low Mach number limits for the Navier-Stokes-Fourier- $P_1$ approximation radiation model in $\mathrm{T}^3$, which is derived under the two physical approximations, that is, "gray approximation" such that the absorbing and scattering coefficients are independent of the frequency, and " $P_1$ approximation" such that the distribution of photons is almost isotropic. We should point out that the results in [1,2,3] are only restricted to the case that the temperature has a small variation, i.e., $\theta = \bar{\theta} + O(\varepsilon)$. Recently, Jiang-Ju-Liao [11] proved the uniform local existence of smooth solutions and the convergence of the model to the incompressible nonhomogeneous Euler system coupled with a diffusion equation, by nonequilibrium-diffusion limit at low Mach number for the Euler-Fourier- $P_1$ approximation radiation model in $\mathrm{T}^3$, when the effect of the large temperature variation upon the limit is taken into account.

For the analysis that we have in mind, it is convenient to introduce the nondimensional form of the system (1.1). For that purpose, we consider the following nondimensionalization:

$$\hat{t} = \frac{t}{T_\infty},\ \mathbf{x} = \frac{\mathbf{x}}{L_\infty},\ \hat{\rho} = \frac{\rho}{\rho_\infty},\ \mathbf{u} = \frac{\mathbf{u}}{U_\infty},\ \hat{P} = \frac{P}{P_\infty},\ \hat{\mu} = \frac{\mu}{\mu_\infty},\ \hat{\lambda} = \frac{\lambda}{\mu_\infty},$$

$$\hat{\theta} = \frac{\theta}{\theta_\infty},\ \hat{e} = \frac{e}{e_\infty},\ \hat{\kappa} = \frac{\kappa}{\kappa_\infty},\ \hat{n} = \frac{n}{n_\infty},\ \sigma_a = \frac{\sigma_a}{\sigma_{a,\infty}},\ \hat{\tilde{\sigma}} = \frac{\tilde{\sigma}}{\tilde{\sigma}_\infty},\ \hat{\nu} = \frac{\nu}{L_\infty},$$

where $T_\infty, L_\infty, \rho_\infty, U_\infty, P_\infty, \mu_\infty, \theta_\infty, e_\infty$ and $\kappa_\infty$ denote the reference hydrodynamical quantities (time, length, density, mean velocity, viscosity, temperature, pressure, energy, heat-conductivity), and $n_\infty, \sigma_{a,\infty}, \tilde{\sigma}_\infty$ represent the reference radiative quantities (radiative intensity, absorption and scaled absorption coefficients). Moreover, we also assume the compatibility conditions $P_\infty = \rho_\infty e_\infty, \kappa_\infty = L_\infty P_\infty U_\infty$, and we denote by

$$S_r := \frac{L_\infty}{T_\infty U_\infty}, M_a := \frac{U_\infty}{\sqrt{P_\infty / \rho_\infty}}, Re := \frac{U_\infty \rho_\infty L_\infty}{\mu_\infty}, \mathrm{C} := \frac{c}{U_\infty},$$

the Strouhal, Mach, Reynolds, and ``infrarelativistic" numbers corresponding to hydrodynamics, and by

$$\mathrm{L}_{\tilde{\sigma}} := \frac{\tilde{\sigma}_\infty \theta_\infty^4 L_\infty}{P_\infty U_\infty},\ \mathrm{L}_{\sigma_a} := \frac{\sigma_{a,\infty} n_\infty L_\infty}{P_\infty U_\infty},\ \mathrm{P} := \frac{\rho_\infty}{n_\infty},$$

the various parameters corresponding to radiation. Using the above scaling, we can rewrite (1.1) for the renormalized unknowns (still denoted by as $(\rho, \mathbf{u}, \theta, n, P, e)$) follows :

$$\begin{aligned}
&S_r \rho_t + \mathrm{div}(\rho \mathbf{u}) = 0,\\
&S_r \rho \mathbf{u}_t + \rho(\mathbf{u}\cdot\nabla)\mathbf{u} + \frac{1}{M_a^2}\nabla P = \frac{1}{Re}\big(\mu\Delta\mathbf{u} + (\mu+\lambda)\nabla \mathrm{div}\mathbf{u}\big),\\
&S_r \rho e_t + \rho\mathbf{u}\cdot\nabla e + P\mathrm{div}\mathbf{u}\\
&\quad = \kappa\Delta\theta + M_a^2 Re\big(2\mu\mathbf{D}(\mathbf{u}):\mathbf{D}(\mathbf{u}) + \lambda(\mathrm{div}\mathbf{u})^2\big) - \mathrm{L}_{\tilde{\sigma}}\tilde{\sigma}\theta^4 + \mathrm{L}_{\sigma_a}\sigma_a n,\\
&\frac{S_r}{\mathrm{C}} n_t - \nu\Delta n = \mathrm{P}\,\mathrm{L}_{\tilde{\sigma}}\tilde{\sigma}\theta^4 - \mathrm{P}\,\mathrm{L}_{\sigma_a}\sigma_a n.
\end{aligned} \tag{1.2}$$

The main purpose of this paper is to investigate the low Mach number limit of the system (1.2) at equilibrium-diffusion regime. Therefore, we use $M_a = \epsilon$ and $\mathrm{L}_{\tilde{\sigma}} = \mathrm{L}_{\sigma_a} = 1$ because the scattering effects vanish in diffusion approximation models. We also assume that the flow is strongly underrelativistic, i.e., $\mathrm{C} = \epsilon^{-1}$. Furthermore, for the sake of simplicity, we put $S_r = Re = \mathrm{P} = 1$. Then, we obtain from (1.2) that

$$
\begin{aligned}
&\rho_t^\varepsilon + \operatorname{div}(\rho^\varepsilon \mathbf{u}^\varepsilon) = 0, \\
&\rho^\varepsilon \mathbf{u}_t^\varepsilon + \rho^\varepsilon (\mathbf{u}^\varepsilon \cdot \nabla)\mathbf{u}^\varepsilon + \frac{1}{\varepsilon^2}\nabla P^\varepsilon = \mu \Delta \mathbf{u}^\varepsilon + (\mu+\lambda)\nabla \operatorname{div}\mathbf{u}^\varepsilon, \\
&\rho^\varepsilon e_t^\varepsilon + \rho^\varepsilon \mathbf{u}^\varepsilon \cdot \nabla e^\varepsilon + P^\varepsilon \operatorname{div}\mathbf{u}^\varepsilon \\
&\qquad = \kappa \Delta \theta^\varepsilon + \varepsilon^2 \big(2\mu \mathbf{D}(\mathbf{u}^\varepsilon) : \mathbf{D}(\mathbf{u}^\varepsilon) + \lambda (\operatorname{div}\mathbf{u}^\varepsilon)^2\big) - \tilde{\sigma}(\theta^\varepsilon)^4 + \sigma_a n^\varepsilon, \\
&\varepsilon n_t^\varepsilon - \nu \Delta n^\varepsilon = \tilde{\sigma}(\theta^\varepsilon)^4 - \sigma_a n^\varepsilon.
\end{aligned} \tag{1.3}
$$

We consider the low Mach number limit for the diffusion approximation model (1.3) in 3-D with the initial data

$$
(\rho^\varepsilon, \mathbf{u}^\varepsilon, \theta^\varepsilon, n^\varepsilon)(\mathbf{x}, 0) = (\rho_0^\varepsilon, \mathbf{u}_0^\varepsilon, \theta_0^\varepsilon, n_0^\varepsilon)(\mathbf{x}) \to (\bar{\rho}, 0, \bar{\theta}, \bar{n}) \quad \text{as } |\mathbf{x}| \to \infty, \tag{1.4}
$$

where $\bar{\rho}, \bar{\theta}, \bar{n}$ are the given positive constants, and with the gas law satisfying

$$
P_\rho(\rho,\theta) > 0, \quad e_\theta(\rho,\theta) > 0 \quad \text{for any } \rho > 0, \theta > 0, \tag{1.5}
$$

which is natural and more general than ideal polytropic gas

$$
P = R\rho\theta, \qquad e = c_v \theta,
$$

where the parameters $R > 0$ and $c_v > 0$ are the gas constant and the heat capacity at constant volume, respectively.

Now we state the main results of this paper.

**Theorem 1.1** (***Local solutions and low Mach number limite***) *Assume that* (1.5) *and*

$$
\begin{aligned}
&(\rho_0^\varepsilon - \bar{\rho}, \mathbf{u}_0^\varepsilon, \theta_0^\varepsilon - \bar{\theta}, n_0^\varepsilon - \bar{n}) \in H^N(\mathbb{R}^3), \\
&\inf_{\mathbf{x} \in \mathbb{R}^3} \rho_0^\varepsilon(\mathbf{x}) > 0, \quad \operatorname{div}\mathbf{u}_0 = 0, \quad \inf_{\mathbf{x} \in \mathbb{R}^3} \theta_0^\varepsilon(\mathbf{x}) > 0
\end{aligned} \tag{1.6}
$$

*for an integer* $N \geq 3$, *where* $(\bar{n}, \bar{\theta})$ *satisfy the compatibility condition*

$$
\sigma_a \bar{n} = \tilde{\sigma} \bar{\theta}^4. \tag{1.7}
$$

*Moreover, for a positive constant* $M_0$ *independent of* $\varepsilon$, *the functions* $(\rho_0^\varepsilon, \mathbf{u}_0^\varepsilon, \theta_0^\varepsilon, n_0^\varepsilon)$ *are assumed to satisfy*

$$
\left\| \rho_0^\varepsilon \mathbf{u}_0^\varepsilon \right\|_{H^N} + \varepsilon^{-1} \left\| (\rho_0^\varepsilon - \bar{\rho}, \theta_0^\varepsilon - \bar{\theta}) \right\|_{H^N} + \varepsilon^{-\frac{1}{2}} \left\| n_0^\varepsilon - \bar{n} \right\|_{H^N} \leq M_0. \tag{1.8}
$$

*Then the following statements hold.*

*Local solutions and uniform estimates: There exist constants* $T_0$ *and* $C$ *depending on* $M_0$ *but independent of* $\varepsilon$, *such that the unique smooth solution* $(\rho^\varepsilon, \mathbf{u}^\varepsilon, \theta^\varepsilon, n^\varepsilon)$ *of the Cauchy problem* (1.3), (1.4) *exists for all small* $\varepsilon$ *on the time interval* $[0, T_0]$ *with properties:*

$$
\begin{aligned}
&(\rho^\varepsilon - \bar{\rho}, \mathbf{u}^\varepsilon, \theta^\varepsilon - \bar{\theta}, n^\varepsilon - \bar{n}) \in C([0,T_0]; H^N(\mathbb{R}^3)), \\
&\nabla \rho^\varepsilon \in L^2(0,T_0; H^N(\mathbb{R}^3)), \quad \big(\nabla(\rho^\varepsilon \mathbf{u}^\varepsilon), \nabla\theta^\varepsilon, \nabla n^\varepsilon\big) \in L^2(0,T_0; H^N(\mathbb{R}^3))
\end{aligned} \tag{1.9}
$$

*and*

$$
\begin{aligned}
&\left\| (\rho^\varepsilon \mathbf{u}^\varepsilon)(t) \right\|_{H^N}^2 + \varepsilon^{-2} \left\| \big(\rho^\varepsilon - \bar{\rho}, \theta^\varepsilon - \bar{\theta}\big)(t) \right\|_{H^N}^2 + \varepsilon^{-1} \left\| \big(n^\varepsilon - \bar{n}\big)(t) \right\|_{H^N}^2 \\
&+ \int_0^t \left\| \nabla(\rho^\varepsilon \mathbf{u}^\varepsilon) \right\|_{H^N}^2 d\tau + \varepsilon^{-2} \int_0^t \left\| \nabla \rho^\varepsilon \right\|_{H^{N-1}}^2 d\tau + \varepsilon^{-2} \int_0^t \left\| \big(\nabla\theta^\varepsilon, \nabla n^\varepsilon\big) \right\|_{H^N}^2 d\tau \leq C
\end{aligned} \tag{1.10}
$$

*for all* $t\in[0,T_0]$.

*Low Mach number limit: There exists a function* $\mathbf{u}\in L^\infty(0,T_0;H^N(\mathbb{R}^3))\cap L^2(0,T_0;\,H^{N+1}(\mathbb{R}^3))$ *such that*

$$\begin{aligned}&\mathbf{u}^\epsilon\to\mathbf{u}\,(\text{*-weak) in}\,L^\infty(0,T_0;H^N(\mathbb{R}^3)),\\&\nabla\mathbf{u}^\epsilon\to\nabla\mathbf{u}\,(\text{weak) in}\,L^2(0,T_0;H^N(\mathbb{R}^3))\end{aligned}\tag{1.11}$$

*as* $\epsilon\to 0$, *and the function pair* $(\mathbf{u},P)$ *for some* $P\in L^\infty(0,T_0;H^{N-1}(\mathbb{R}^3))\cap L^2(0,T_0;H^N(\mathbb{R}^3))$ *is the unique smooth local solution on the time interval* $[0,T_0]$ *of the incompressible Navier-Stokes equation*

$$\begin{aligned}&\operatorname{div}\mathbf{u}=0,\\&\mathbf{u}_t+(\mathbf{u}\cdot\nabla)\mathbf{u}+\frac{1}{\bar\rho}\nabla P=\frac{\mu}{\bar\rho}\Delta\mathbf{u},\end{aligned}\tag{1.12}$$

*with the initial data*

$$\mathbf{u}(x,0)=\mathbf{u}_0(x),\tag{1.13}$$

*where* $\mathbf{u}_0\in H^N(\mathbb{R}^3)$ *is a function such that*

$$\mathbf{u}_0^\epsilon\to\mathbf{u}_0,(\text{weak) in};H^N(\mathbb{R}^3)\quad\text{as}\quad\epsilon\to 0.$$

**Theorem 1.2 (*Global solutions and low Mach number limit*)** *Assume that* (1.5), (1.6) *and* (1.7). *Moreover, for small positive constant* $\delta_0>0$ *independent of* $\epsilon$, *the functions* $(\rho_0^\epsilon,\mathbf{u}_0^\epsilon,\theta_0^\epsilon,n_0^\epsilon)$ *are assumed to satisfy*

$$\|\mathbf{u}_0^\epsilon\|_{H^3}+\epsilon^{-1}\|(\rho_0^\epsilon-\bar\rho,\theta_0^\epsilon-\bar\theta)\|_{H^3}+\epsilon^{-\frac12}\|n_0^\epsilon-\bar n\|_{H^3}\le\delta_0.\tag{1.14}$$

*Then the following statements hold.*

*Global solutions and uniform estimate: There exist constant* $C$ *independent of* $\epsilon$ *and* $\delta_0$ *such that the Cauchy problem* (1.3), (1.4) *admits a unique solution* $(\rho^\epsilon,\mathbf{u}^\epsilon,\theta^\epsilon,n^\epsilon)$ *on* $[0,\infty)$ *satisfying*

$$\begin{aligned}&(\rho^\epsilon-\bar\rho,\mathbf{u}^\epsilon,\theta^\epsilon-\bar\theta,n^\epsilon-\bar n)\in C([0,\infty),H^N(\mathbb{R}^3)),\\&\nabla\rho^\epsilon\in L^2(0,\infty;H^{N-1}(\mathbb{R}^3)),(\nabla\mathbf{u}^\epsilon,\nabla\theta^\epsilon,\nabla n^\epsilon)\in L^2(0,\infty;H^N(\mathbb{R}^3))\end{aligned}$$

*and*

$$\begin{aligned}&\sup_{t\in\mathbb{R}_+}\|\mathbf{u}^\epsilon(t)\|_{H^N}^2+\epsilon^{-2}\sup_{t\in\mathbb{R}_+}\|(\rho^\epsilon-\bar\rho,\theta^\epsilon-\bar\theta)(t)\|_{H^N}^2+\epsilon^{-1}\sup_{t\in\mathbb{R}_+}\|(n^\epsilon-\bar n)(t)\|_{H^N}^2\\&+\epsilon^{-2}\int_0^\infty\|\nabla\rho^\epsilon\|_{H^{N-1}}^2d\tau+\int_0^\infty\|\nabla\mathbf{u}^\epsilon\|_{H^N}^2d\tau+\epsilon^{-2}\int_0^\infty\|(\nabla\theta^\epsilon,\nabla n^\epsilon)\|_{H^N}^2d\tau\\&\le C(\|\mathbf{u}_0^\epsilon\|_{H^N}^2+\epsilon^{-2}\|(\rho_0^\epsilon-\bar\rho,\theta_0^\epsilon-\bar\theta)\|_{H^N}^2+\epsilon^{-1}\|n_0^\epsilon-\bar n\|_{H^N}^2).\end{aligned}\tag{1.15}$$

Low Mach number limit: There exists a function $\mathbf{u}\in L^\infty(0,\infty;H^N(\mathbb{R}^3)),\nabla\mathbf{u}\in L^2\,(0,\infty;H^N(\mathbb{R}^3))$ such that

$$\begin{aligned}&\mathbf{u}^\epsilon\to\mathbf{u}\,(\text{*-weak) in}\,L^\infty(0,\infty;H^N(\mathbb{R}^3)),\\&\nabla\mathbf{u}^\epsilon\to\nabla\mathbf{u}\,(\text{weak) in}\,L^2(0,\infty;H^N(\mathbb{R}^3))\end{aligned}$$

as $\epsilon \to 0$, and the function pair $(\mathbf{u}, P)$ for some $P \in L^\infty(0,\infty; H^{N-1}(\mathbb{R}^3)) \cap L^2(0,\infty;\ H^N(\mathbb{R}^3))$ is the unique smooth global solution on $[0,\infty)$ to the Cauchy problem (1.12), (1.13).

**Remark 1.1** *For the global existence of solutions, we only assume that the $H^3$-norm of initial perturbations is small, while the higher-order Sobolev norms can be arbitrarily large* (*see* (1.14)).

**Highlight of this paper.** We briefly review the key analytical techniques. To prove Theorem 1.1 and Theorem 1.2, we will use an elementary energy method, where the essential step is a priori estimates independent of Mach number parameter $\epsilon$ for the perturbative solutions. The main difficulty comes from the strong nonlinear term $(\theta^\epsilon)^4$ due to the radiation effect in (1.1)3,4, while the nonlinear coupling between the Navier-Stokes equations (1.1)1-3 and the radiative transport equation (1.1)4 also bring trouble to close a priori estimates. In order to prove Theorem 1.1, we first reformulate the system (1.3) into (2.3) for the perturbation $(\tilde{n}, \mathbf{m}, \zeta, \vartheta)$ of the density $\rho^\epsilon$, the momentum $\rho^\epsilon \mathbf{u}^\epsilon$, the temperature $\theta^\epsilon$ and radiation field $n^\epsilon$. Here, a key is to add a leading term $(4\tilde{\sigma}\bar{\theta}^3\zeta - \sigma_a \vartheta)$ which is representing an interaction between the perturbations of the temperature $\theta^\epsilon = \bar{\theta} + \zeta$ and radiation field $n^\epsilon = \bar{n} + \vartheta$ (see (2.2)). The main step in obtaining the existence and low Mach number limit is the estimation (2.12) independent of $\epsilon$ for the linearized system (2.8). Then, we can the desired results by using the linearization of the original problem and the boundness of the iterative sequences of solutions for the reformulated system (Subsection (2.3)). Also, the main steps to prove Theorem 1.2 are the energy estimations (3.8) and (3.24) independent of $\epsilon$, which are possible to close a priori estimate needed for the continuation argument of local solutions and the low Much number limit (Subsection (3.2)).

**Notation.** In this paper, $L^p(\mathbb{R}^3)$ and $W_p^k(\mathbb{R}^3)$ denote the usual Lebesgue and Sobolev spaces on $\mathbb{R}^3$, with norms $\|\cdot\|_{L^p}$ and $\|\cdot\|_{W_p^k}$, respectively. When $p=2$, we denote $W_p^k(\mathbb{R}^3)$ by $H^k(\mathbb{R}^3)$ with the norm $\|\cdot\|_{H^k}$ and $\|\cdot\|_{H^0} = \|\cdot\|$ will be used to denote the usual $L^2-$norm. The notation $\|(A_1, A_2, \cdots, A_l)\|_{H^k}$ means the summation of $\|A_i\|_{H^k}$ from $i=1$ to $i=l$. For an integer $m$, the symbol $\nabla^m$ denotes the summation of all terms $D^\alpha$ with the multi-index $\alpha$ satisfying $|\alpha| = m$. We omit the spatial domain $\mathbb{R}^3$ in integrals for convenience.

Before finishing this section, we will recall the following useful Lemmas which we will use extensively.

**Lemma 1.1** (*Gagliardo-Nirenberg inequality,* [15]) *Let $l, s$ and $k$ be any real numbers satisfying $0 \le l, s \le k$, and let $p, r, q \in [1,\infty]$ and $0 \le \theta \le 1$ such that*

$$\frac{l}{3} - \frac{1}{p} = \left(\frac{s}{3} - \frac{1}{r}\right)(1-\theta) + \left(\frac{k}{3} - \frac{1}{q}\right)\theta.$$

*Then, for any $u \in W_q^k(\mathbb{R}^3)$, we have*

$$\|\nabla^l u\|_{L^p} \le C \|\nabla^s u\|_{L^r}^{1-\theta} \|\nabla^k u\|_{L^q}^{\theta}. \tag{1.16}$$

**Lemma 1.2** (*see* [4]) *Let $\Omega = \mathbb{R}^d$, $s_1 \ge s$ and $s_2 \ge s$ be such that either*

$$s_1 + s_2 - s \ge d\left(\frac{1}{q_1} + \frac{1}{q_2} - \frac{1}{q}\right) \ge 0, \quad c, \quad j = 1, 2$$

*or*

$$s_1+s_2-s>d(\frac{1}{q_1}+\frac{1}{q_2}-\frac{1}{q})\geq 0,\quad s_j-s\geq d(\frac{1}{q_j}-\frac{1}{q}),\ j=1,2,$$

*then* $(u,v)\mapsto u\cdot v$ *is a continuous bilinear map from* $W^{s_1}_{q_1}(\Omega)\times W^{s_1}_{q_1}(\Omega)$ *into* $W^{s}_{q}(\Omega)$\$.

**Lemma 1.3** ([17, *Lemma 2.5*]) *Let* $f(\varphi)$ *and* $f(\varphi,w)$ $f(\varphi,w)$ *be smooth functions of* $\varphi$ *and* $(\varphi,w)$, *respectively, with bounded derivatives of any order, and* $\|\varphi\|_{L^\infty(\mathbb{R}^3)}+\|w\|_{L^\infty(\mathbb{R}^3)}\leq C$. *Then for any integer* $m\geq 1$, *we have*

$$\begin{aligned}&\left\|\nabla^m f(\varphi)\right\|_{L^p}\leq C\left\|\nabla^m\varphi\right\|_{L^p},\\&\left\|\nabla^m f(\varphi,w)\right\|_{L^p}\leq C\left\|\nabla^m(\varphi,w)\right\|_{L^p},\end{aligned}\tag{1.17}$$

*for any* $1\leq p\leq\infty$, *where* $C$ *may depend on* $f$ *and* $m$.

**Lemma 1.4** ( [17, *Lemma* 2.6]) Let $\alpha$ be any multi-index with $|\alpha|=k$ and $1<p<\infty$. *Then there exists some constant* $C>0$ *such that*

$$\begin{aligned}&\left\|D^\alpha(fg)\right\|_{L^p}\leq C(\|f\|_{L^{p_1}}\left\|\nabla^k g\right\|_{L^{p_2}}+\left\|\nabla^k f\right\|_{L^{p_3}}\|g\|_{L^{p_4}}),\\&\left\|[D^\alpha,f]g\right\|_{L^p}\leq C\left\|\nabla f\right\|_{L^{p_1}}\left\|\nabla^{k-1}g\right\|_{L^{p_2}}+\left\|\nabla^k f\right\|_{L^{p_3}}\|\,\|g\|\,\|_{L^{p_4}}),\end{aligned}\tag{1.18}$$

*where* $f,g\in S$ *is the Schwartz class,* $1\leq p_2,p_3\leq\infty$ *such that* $\frac{1}{p}=\frac{1}{p_1}+\frac{1}{p_2}=\frac{1}{p_3}+\frac{1}{p_4}$, *and* $[D^\alpha,f]g=D^\alpha(fg)-fD^\alpha g$.

## 2. Local solution and low Mach number limit

In this section, we prove Theorem 1.1. To this end, we first reformulate the system (1.3) to a new system. For brevity, we denote $(\rho^\epsilon,\mathbf{u}^\epsilon,\theta^\epsilon,n^\epsilon)$, $P^\epsilon=P(\rho^\epsilon,\theta^\epsilon)$ and $e^\epsilon=e(\rho^\epsilon,\theta^\epsilon)$ in (1.3) by $(\rho,\mathbf{u},\theta,n)$, $P=P(\rho,\theta)$ and $e=e(\rho,\theta)$, respectively.

### 2.1 Reformulated system

By using the following thermodynamical relation $-\rho^2 e_\rho(\rho,\theta)=\theta P_\theta(\rho,\theta)-P(\rho,\theta)$ (see [5,(1.6)]), we rewrite $(1.3)_3$ as

$$\begin{aligned}&\theta_t+\mathbf{u}\cdot\nabla\theta+\frac{\theta P_\theta(\rho,\theta)}{\rho e_\theta(\rho,\theta)}\mathrm{div}\mathbf{u}-\frac{\kappa}{\rho e_\theta(\rho,\theta)}\Delta\theta\\&\qquad=\frac{2\mu\epsilon^2\mathbf{D}(\mathbf{u}):\mathbf{D}(\mathbf{u})}{\rho e_\theta(\rho,\theta)}+\frac{\lambda\epsilon^2(\mathrm{div}\mathbf{u})^2}{\rho e_\theta(\rho,\theta)}-\frac{\tilde{\sigma}\theta^4-\sigma_a n}{\rho e_\theta(\rho,\theta)}.\end{aligned}\tag{2.1}$$

Setting

$$\tilde{n}=\frac{\rho-\bar{\rho}}{\bar{\rho}},\quad \mathbf{m}=\frac{\rho\mathbf{u}}{\bar{\rho}}\quad \zeta=\theta-\bar{\theta}\quad\text{and}\quad \vartheta=n-\bar{n},$$

we get

$$\tilde{\sigma}(\zeta+\bar{\theta})^4-\sigma_a(\vartheta+\bar{n})=4\tilde{\sigma}\bar{\theta}^3\zeta-\sigma_a\vartheta+6\tilde{\sigma}\bar{\theta}^2\zeta^2+4\tilde{\sigma}\bar{\theta}\zeta^3+\tilde{\sigma}\zeta^4\tag{2.2}$$

due to (1.7).We also define the scaled viscosity coefficients $\bar{\mu}=\frac{\mu}{\bar{\rho}}$ and $\bar{\lambda}=\frac{\lambda}{\bar{\rho}}$. Assuming that the density $\rho$ is bounded away from zero, we rewrite the diffusion approximation model (1.3) as follows:

$$
\begin{aligned}
&\tilde{n}_t+\mathrm{div}\mathbf{m}=0,\quad \mathbf{x}\in\Box^3, t>0,\\
&\mathbf{m}_t-\bar{\mu}\mathrm{div}(\frac{\nabla\mathbf{m}}{1+\tilde{n}})-(\bar{\lambda}+\bar{\mu})\nabla(\frac{\mathrm{div}\mathbf{m}}{1+\tilde{n}})\\
&\qquad+\frac{\bar{\rho}P_\rho(\bar{\rho},\bar{\theta})}{\grave{o}^2}\nabla\tilde{n}+\frac{P_\theta(\bar{\rho},\bar{\theta})}{\grave{o}^2}\nabla\zeta=\mathbf{G}_1(\tilde{n},\mathbf{m},\zeta),\\
&\zeta_t+\frac{\bar{\theta}P_\theta(\bar{\rho},\bar{\theta})}{\bar{\rho}e_\theta(\bar{\rho},\bar{\theta})}\mathrm{div}\mathbf{m}=\frac{\kappa}{\bar{\rho}e_\theta(\bar{\rho},\bar{\theta})}\Delta\zeta-\frac{4\tilde{\sigma}\bar{\theta}^3\zeta-\sigma_a\vartheta}{\bar{\rho}e_\theta(\bar{\rho},\bar{\theta})}+G_2(\tilde{n},\mathbf{m},\zeta),\\
&\grave{o}\vartheta_t-\nu\Delta\vartheta=4\tilde{\sigma}\bar{\theta}^3\zeta-\sigma_a\vartheta+G_3(\zeta),
\end{aligned}
\tag{2.3}
$$

where we used

$$
1+\tilde{n}=\frac{\rho}{\bar{\rho}},\quad \mathbf{u}=\frac{\mathbf{m}}{1+\tilde{n}},
$$

and

$$
\begin{aligned}
\mathbf{G}_1(\tilde{n},\mathbf{m},\zeta)=&-\mathrm{div}\left(\frac{\mathbf{m}\otimes\mathbf{m}}{1+\tilde{n}}\right)+\bar{\mu}\mathrm{div}\left(\mathbf{m}\otimes\nabla\left(\frac{1}{1+\tilde{n}}\right)\right)\\
&+(\bar{\lambda}+\bar{\mu})\nabla\left(\mathbf{m}\cdot\nabla\left(\frac{1}{1+\tilde{n}}\right)\right)+\frac{\bar{\rho}h_1(\tilde{n},\zeta)}{\grave{o}^2}\nabla\tilde{n}+\frac{h_2(\tilde{n},\zeta)}{\grave{o}^2}\nabla\zeta,
\end{aligned}
\tag{2.4}
$$

$$
h_1(\tilde{n},\zeta)=P_\rho(\bar{\rho},\bar{\theta})-P_\rho(\rho,\theta),\ h_2(\tilde{n},\zeta)=P_\theta(\bar{\rho},\bar{\theta})-P_\theta(\rho,\theta),
$$

$$
\begin{aligned}
G_2(\tilde{n},\mathbf{m},\zeta)=&-\left(\frac{\mathbf{m}}{1+\tilde{n}}\right)\cdot\nabla\zeta-\kappa h_3(\tilde{n},\zeta)\Delta\zeta-h_3(\tilde{n},\zeta)\left(4\tilde{\sigma}\bar{\theta}^3\zeta-\sigma_a\vartheta\right)\\
&+\frac{2\mu\grave{o}^2\mathbf{D}\left(\frac{\mathbf{m}}{1+\tilde{n}}\right):\mathbf{D}\left(\frac{\mathbf{m}}{1+\tilde{n}}\right)}{\rho e_\theta(\rho,\theta)}+\frac{\lambda\grave{o}^2\left(\mathrm{div}\left(\frac{\mathbf{m}}{1+\tilde{n}}\right)\right)^2}{\rho e_\theta(\rho,\theta)}+h_4(\tilde{n},\zeta)\mathrm{div}\left(\frac{\mathbf{m}}{1+\tilde{n}}\right)\\
&-\frac{h_5(\zeta)\zeta}{\rho e_\theta(\rho,\theta)}+\frac{\bar{\theta}P_\theta(\bar{\rho},\bar{\theta})}{e_\theta(\bar{\rho},\bar{\theta})}\left(\frac{\tilde{n}}{1+\tilde{n}}\mathrm{div}\mathbf{m}-\mathbf{m}\cdot\nabla\left(\frac{1}{1+\tilde{n}}\right)\right),
\end{aligned}
\tag{2.5}
$$

$$
h_3(\tilde{n},\zeta)=\frac{1}{\bar{\rho}e_\theta(\bar{\rho},\bar{\theta})}-\frac{1}{\rho e_\theta(\rho,\theta)},\ h_4(\tilde{n},\zeta)=\frac{\bar{\theta}P_\theta(\bar{\rho},\bar{\theta})}{\bar{\rho}e_\theta(\bar{\rho},\bar{\theta})}-\frac{\theta P_\theta(\rho,\theta)}{\rho e_\theta(\rho,\theta)},
$$

$$
h_5(\zeta)=6\tilde{\sigma}\bar{\theta}^2\zeta+4\tilde{\sigma}\bar{\theta}\zeta^2+\tilde{\sigma}\zeta^3
$$

$$
G_3(\zeta)=h_5(\zeta)\zeta. \tag{2.6}
$$

Also, the initial condition (1.4) is reformulated into

$$(\varrho,\mathbf{m},\zeta,\vartheta)(\mathbf{x},0)=(\varrho_0,\mathbf{m}_0,\zeta_0,\vartheta_0)(\mathbf{x})$$
$$\equiv\left(\frac{\rho_0^{\epsilon}(\mathbf{x})-\bar{\rho}}{\bar{\rho}},\frac{\rho_0^{\epsilon}(\mathbf{x})\mathbf{u}_0^{\epsilon}(\mathbf{x})}{\bar{\rho}},\theta_0^{\epsilon}(\mathbf{x})-\bar{\theta},n_0^{\epsilon}(\mathbf{x})-\bar{n}\right). \tag{2.7}$$

The proof of the local existence theorem for the reformulated system (2.3), (2.7) relies upon the study of a linearized system and the contraction mapping principle.

## 2.2 Estimates for the linearized system

In this subsection, we study the following linearized system

$$\begin{cases}\varrho_t+\mathrm{div}\mathbf{m}=0, \qquad \mathbf{x}\in\mathbb{R}^3,t>0,\\ \mathbf{m}_t-\bar{\mu}\mathrm{div}\left(A\nabla\mathbf{m}\right)-\left(\bar{\lambda}+\bar{\mu}\right)\nabla\left(A\mathrm{div}\mathbf{m}\right)+\dfrac{\bar{\rho}P_{\rho}\left(\bar{\rho},\bar{\theta}\right)}{\epsilon^2}\nabla\varrho+\dfrac{P_{\theta}\left(\bar{\rho},\bar{\theta}\right)}{\epsilon^2}\nabla\zeta=\mathbf{G}_1,\\ \zeta_t+\dfrac{\bar{\theta}P_{\theta}\left(\bar{\rho},\bar{\theta}\right)}{\bar{\rho}e_{\theta}\left(\bar{\rho},\bar{\theta}\right)}\mathrm{div}\mathbf{m}=\dfrac{\kappa}{\bar{\rho}e_{\theta}\left(\bar{\rho},\bar{\theta}\right)}\Delta\zeta-\dfrac{4\tilde{\sigma}\bar{\theta}^3\zeta-\sigma_a\vartheta}{\bar{\rho}e_{\theta}(\bar{\rho},\bar{\theta})}+G_2,\\ \epsilon\vartheta_t-\nu\Delta\vartheta=4\tilde{\sigma}\bar{\theta}^3\zeta-\sigma_a\vartheta+G_3,\\ (\varrho,\mathbf{m},\zeta,\vartheta)(\mathbf{x},0)=\left(\varrho_0,\mathbf{m}_0,\zeta_0,\vartheta_0\right)(\mathbf{x}), \quad \mathbf{x}\in\Omega,\end{cases} \tag{2.8}$$

where $A=A(\mathbf{x},t)$ is the scalar function satisfying

$$0<M_0\le A(\mathbf{x},t)\le M_1<\infty, \tag{2.9}$$

where $M_0,M_1$ are positive constants independently of $(\mathbf{x},t)$.

We first will consider an estimate for the linearized system (2.8).

**Theorem 2.1** *Let* $0<T<\infty$. *Assume that* $(\varrho_0,\mathbf{m}_0,\zeta_0,\vartheta_0)\in H^N(\mathbb{R}^3)$ *for an integer* $N\ge 3$. *Also, suppose that* (2.9) *and*

$$A\in L^{\infty}\left(0,T;H^N\left(\mathbb{R}^3\right)\right), \quad \mathbf{G}_1,G_2,G_3\in L^2\left(0,T;H^{N-1}\left(\mathbb{R}^3\right)\right). \tag{2.10}$$

*Let*

$$\begin{aligned}&\left(\varrho,\mathbf{m},\zeta,\vartheta\right)\in C\left(\left[0,T\right];H^N\left(\mathbb{R}^3\right)\right),\\ &\nabla\varrho\in L^2\left(0,T;H^{N-1}\left(\mathbb{R}^3\right)\right), \quad \left(\nabla\mathbf{m},\nabla\zeta,\nabla\vartheta\right)\in L^2\left(0,T;H^N\left(\mathbb{R}^3\right)\right)\end{aligned} \tag{2.11}$$

*and* $(\varrho,\mathbf{m},\zeta,\vartheta)$ *be a solution of the system* (2.8). *Then there exist positive constants* $C_0$ *and* $c_0$, *independent of* $\epsilon$ *and* $t$, *such that the following inequality holds:*

$$
\begin{aligned}
&\left\|\left(\mathbf{m},\grave{o}^{-1}\tilde{n},\grave{o}^{-1}\zeta,\grave{o}^{-\frac{1}{2}}\vartheta\right)(t)\right\|_{H^N}^2+\grave{o}^{-2}\int_0^t\|\nabla\tilde{n}\|_{H^{N-1}}^2 d\tau\\
&+\int_0^t\|\nabla\mathbf{m}\|_{H^N}^2 d\tau+\grave{o}^{-2}\int_0^t\left\|\left(\nabla\zeta,\nabla\vartheta,4\tilde{\sigma}\bar{\theta}^3\zeta-\sigma_a\vartheta\right)\right\|_{H^N}^2 d\tau\\
&\le C_0\left[\left\|(\mathbf{m},\grave{o}^{-1}\tilde{n},\grave{o}^{-1}\zeta,\grave{o}^{-\frac{1}{2}}\vartheta)(0)\right\|_{H^N}^2+\int_0^t\left\|\left(\mathbf{G}_1,\grave{o}^{-1}G_2,\grave{o}^{-1}G_3\right)\right\|_{H^{N-1}}^2 d\tau\right]\times\\
&\times\left[1+e^{C_0\int_0^T\left(1+\|A\|_{H^N}^2\right)d\tau}\int_0^T\left(1+\|A\|_{H^N}^2\right)d\tau\right]
\end{aligned}
\tag{2.12}
$$

*for any* $t\in[0,T]$.

**Proof**. Applying $\nabla^k$ to (2.8) yields

$$
\begin{aligned}
&\nabla^k\tilde{n}_t+\mathrm{div}\nabla^k\mathbf{m}=0,\\
&\nabla^k\mathbf{m}_t-\bar{\mu}\mathrm{div}\nabla^k(A\nabla\mathbf{m})-(\bar{\lambda}+\bar{\mu})\nabla\nabla^k(A\mathrm{div}\mathbf{m})\\
&\qquad+\frac{\bar{\rho}P_\rho\left(\bar{\rho},\bar{\theta}\right)}{\grave{o}^2}\nabla^{k+1}\tilde{n}+\frac{P_\theta\left(\bar{\rho},\bar{\theta}\right)}{\grave{o}^2}\nabla\zeta^{k+1}=\nabla^k\mathbf{G}_1,\\
&\nabla^k\zeta_t+\frac{\bar{\theta}P_\theta\left(\bar{\rho},\bar{\theta}\right)}{\bar{\rho}e_\theta\left(\bar{\rho},\bar{\theta}\right)}\mathrm{div}\nabla^k\mathbf{m}=\frac{\kappa}{\bar{\rho}e_\theta\left(\bar{\rho},\bar{\theta}\right)}\Delta\nabla^k\zeta-\frac{4\tilde{\sigma}\bar{\theta}^3\nabla^k\zeta-\sigma_a\nabla^k\vartheta}{\bar{\rho}e_\theta\left(\bar{\rho},\bar{\theta}\right)}+\nabla^k G_2,\\
&\grave{o}\nabla^k\vartheta_t-\nu\Delta\nabla^k\vartheta=4\tilde{\sigma}\bar{\theta}^3\nabla^k\zeta-\sigma_a\nabla^k\vartheta+\nabla^k G_3,
\end{aligned}
\tag{2.13}
$$

where $k=0,1,\cdots,N$.

Multiplying $(2.13)_2$ by $\nabla^k\mathbf{m}$ and using $(2.13)_1$, we have

$$
\begin{aligned}
&\frac{1}{2}\frac{d}{dt}\left(\left\|\nabla^k\mathbf{m}\right\|^2+\frac{\bar{\rho}P_\rho(\bar{\rho},\bar{\theta})}{\grave{o}^2}\left\|\nabla^k\tilde{n}\right\|^2\right)+\int A\left(\bar{\mu}\left|\nabla^{k+1}\mathbf{m}\right|^2+\left(\bar{\lambda}+\bar{\mu}\right)\left|\mathrm{div}\nabla^k\mathbf{m}\right|^2\right)d\mathbf{x}\\
&=\frac{P_\theta\left(\bar{\rho},\bar{\theta}\right)}{\grave{o}^2}\int\nabla^k\zeta\,\mathrm{div}\nabla^k\mathbf{m}d\mathbf{x}+\int\nabla^k\mathbf{G_1}\cdot\nabla^k\mathbf{m}d\mathbf{x}\\
&\quad+\bar{\mu}\int\left[\nabla^k,A\right]\nabla\mathbf{m}:\nabla^{k+1}\mathbf{m}d\mathbf{x}+\left(\bar{\lambda}+\bar{\mu}\right)\int\left[\nabla^k,A\right]\mathrm{div}\mathbf{m}\,\mathrm{div}\nabla^k\mathbf{m}d\mathbf{x},
\end{aligned}
\tag{2.14}
$$

where we used

$$
\begin{aligned}
&\nabla^k\left(A\nabla\mathbf{m}\right)=A\nabla\nabla^k\mathbf{m}-\left[\nabla^k,A\right]\nabla\mathbf{m},\\
&\nabla^k(A\mathrm{div}\mathbf{m})=A\mathrm{div}\nabla^k\mathbf{m}-\left[\nabla^k,A\right]\mathrm{div}\mathbf{m}.
\end{aligned}
$$

Multiplying $(2.13)_3$ by $\frac{\bar{\rho}e_\theta(\bar{\rho},\bar{\theta})}{\grave{o}^2\bar{\theta}}\nabla^k\zeta$, we have

$$
\begin{aligned}
&\frac{\bar{\rho}e_\theta(\bar{\rho},\bar{\theta})}{2\grave{o}^2\bar{\theta}}\frac{d}{dt}\left\|\nabla^k\zeta\right\|^2+\frac{\kappa}{\grave{o}^2\bar{\theta}}\left\|\nabla^{k+1}\zeta\right\|^2+\frac{1}{\grave{o}^2\bar{\theta}}\int\left(4\tilde{\sigma}\bar{\theta}^3\nabla^k\zeta-\sigma_a\nabla^k\vartheta\right)\nabla^k\zeta d\mathbf{x}\\
&\qquad=-\frac{P_\theta(\bar{\rho},\bar{\theta})}{\grave{o}^2}\int\nabla^k\zeta\,\mathrm{div}\nabla^k\mathbf{m}d\mathbf{x}+\frac{\bar{\rho}e_\theta(\bar{\rho},\bar{\theta})}{\grave{o}^2\bar{\theta}}\int\nabla^k G_2\nabla^k\zeta d\mathbf{x}.
\end{aligned}
\tag{2.15}
$$

Multiplying $(2.13)_4$ by $\frac{\sigma_a}{4\text{ò}^2\tilde{\sigma}\bar{\theta}^4}\nabla^k\vartheta$, we have

$$\frac{\sigma_a}{8\text{ò}\tilde{\sigma}\bar{\theta}^4}\frac{d}{dt}\left\|\nabla^k\vartheta\right\|^2+\frac{\sigma_a\eta}{4\text{ò}^2\tilde{\sigma}\bar{\theta}^4}\left\|\nabla^{k+1}\vartheta\right\|^2-\frac{\sigma_a}{4\text{ò}^2\tilde{\sigma}\bar{\theta}^4}\int\left(4\tilde{\sigma}\bar{\theta}^3\nabla^k\zeta-\sigma_a\nabla^k\vartheta\right)\nabla^k\vartheta d\mathbf{x}$$
$$=\frac{\sigma_a}{4\text{ò}^2\tilde{\sigma}\bar{\theta}^4}\int\nabla^k G_3\nabla^k\vartheta d\mathbf{x}. \tag{2.16}$$

Adding (2.14)-(2.16) yields that

$$\frac{1}{2}\frac{d}{dt}\left(\|\nabla^k\mathbf{m}\|^2+\frac{\bar{\rho}P_\rho(\bar{\rho},\bar{\theta})}{\text{ò}^2}\|\nabla^k\text{ñ}\|^2+\frac{\bar{\rho}e_\theta(\bar{\rho},\bar{\theta})}{\text{ò}^2\bar{\theta}}\|\nabla^k\zeta\|^2+\frac{\sigma_a}{4\text{ò}\tilde{\sigma}\bar{\theta}^4}\|\nabla^k\vartheta\|^2\right)$$
$$+\int A\left(\bar{\mu}|\nabla^{k+1}\mathbf{m}|^2+(\bar{\lambda}+\bar{\mu})|\mathrm{div}\nabla^k\mathbf{m}|^2\right)d\mathbf{x}+\frac{\kappa}{\text{ò}^2\bar{\theta}}\|\nabla^{k+1}\zeta\|^2 \tag{2.17}$$
$$+\frac{\sigma_a\eta}{4\text{ò}^2\tilde{\sigma}\bar{\theta}^4}\|\nabla^{k+1}\vartheta\|^2+\frac{\sigma_a}{4\text{ò}^2\tilde{\sigma}\bar{\theta}^4}\|(4\tilde{\sigma}\bar{\theta}^3\nabla^k\zeta-\sigma_a\nabla^k\vartheta)\|^2=J_k^1+J_k^2,$$

where

$$J_k^1=\bar{\mu}\int\left[\nabla^k,A\right]\nabla\mathbf{m}:\nabla^{k+1}\mathbf{m}d\mathbf{x}+\left(\bar{\lambda}+\bar{\mu}\right)\int\left[\nabla^k,A\right]\mathrm{div}\mathbf{m}\,\mathrm{div}\nabla^k\mathbf{m}d\mathbf{x},$$
$$J_k^2=\int\nabla^k\mathbf{G_1}\cdot\nabla^k\mathbf{m}d\mathbf{x}+\frac{\bar{\rho}e_\theta(\bar{\rho},\bar{\theta})}{\text{ò}^2\bar{\theta}}\int\nabla^k G_2\nabla^k\zeta d\mathbf{x}+\frac{\sigma_a}{4\text{ò}^2\tilde{\sigma}\bar{\theta}^4}\int\nabla^k G_3\nabla^k\vartheta d\mathbf{x}. \tag{2.18}$$

Summing up (2.17) for $k=0,1,\cdots,N$ and using (2.9), we get

$$\frac{1}{2}\frac{d}{dt}\left(\|\mathbf{m}\|_{H^N}^2+\frac{\bar{\rho}P_\rho(\bar{\rho},\bar{\theta})}{\text{ò}^2}\|\text{ñ}\|_{H^N}^2+\frac{\bar{\rho}e_\theta(\bar{\rho},\bar{\theta})}{\text{ò}^2\bar{\theta}}\|\zeta\|_{H^N}^2+\sigma_a 4\text{ò}\tilde{\sigma}\bar{\theta}^4\|\vartheta\|_{H^N}^2\right)$$
$$+c_0\|\nabla\mathbf{m}\|_{H^N}^2+\frac{\kappa}{\text{ò}^2\bar{\theta}}\|\nabla\zeta\|_{H^N}^2+\frac{\sigma_a\eta}{4\text{ò}^2\tilde{\sigma}\bar{\theta}^4}\|\nabla\vartheta\|_{H^N}^2 \tag{2.19}$$
$$+\frac{\sigma_a}{4\text{ò}^2\tilde{\sigma}\bar{\theta}^4}\left\|\left(4\tilde{\sigma}\bar{\theta}^3\zeta-\sigma_a\vartheta\right)\right\|_{H^N}^2\le\sum_{k=0}^N J_k^1+\sum_{k=0}^N J_k^2.$$

Using Holder inequality and (1.18), we obtain from (2.18) that

$$\sum_{k=0}^N J_k^1\le C\sum_{k=0}^N\left(\|\nabla A\|_{L^\infty}\|\nabla^k\mathbf{m}\|+\|\nabla\mathbf{m}\|_{L^\infty}\|\nabla^k A\|\right)\|\nabla^{k+1}\mathbf{m}\|$$
$$+C\sum_{k=0}^N\left(\|\nabla A\|_{L^\infty}\|\nabla^{k-1}\mathrm{div}\mathbf{m}\|+\|\mathrm{div}\mathbf{m}\|_{L^\infty}\|\nabla^k A\|\right)\|\nabla^k\mathrm{div}\mathbf{m}\| \tag{2.20}$$
$$\le C\|A\|_{H^N}\|\mathbf{m}\|_{H^N}\|\nabla\mathbf{m}\|_{H^N}.$$

and

$$\sum_{k=0}^N J_k^2\le C\sum_{k=0}^{N-1}\left(\|\nabla^k\mathbf{G}_1\|\|\nabla^k\mathbf{m}\|+\text{ò}^{-2}\|\nabla^k G_2\|\|\nabla^k\zeta\|+\text{ò}^{-2}\|\nabla^k G_3\|\|\nabla^k\vartheta\|\right)$$
$$+\left|\int\nabla^N\mathbf{G}_1\cdot\nabla^{N+1}\mathbf{m}d\mathbf{x}\right|+\text{ò}^{-2}\left|\int\left(\nabla^N G_1\nabla^{N+1}\zeta+\nabla^N G_3\nabla^{N+1}\vartheta\right)d\mathbf{x}\right| \tag{2.21}$$
$$\le C\left(\|\mathbf{G}_1\|_{H^{N-1}}\|\nabla\mathbf{m}\|_{H^N}+\text{ò}^{-2}\|G_2\|_{H^{N-1}}\|\nabla\zeta\|_{H^N}+\text{ò}^{-2}\|G_3\|_{H^{N-1}}\|\nabla\vartheta\|_{H^N}\right).$$

Substituting (2.20)-(2.21) into (2.19), we get

$$
\frac{d}{dt}\left(\|\mathbf{m}\|_{H^N}^2+\frac{\bar{\rho}P_\rho(\bar{\rho},\bar{\theta})}{\grave{o}^2}\|\tilde{n}\|_{H^N}^2+\frac{\bar{\rho}e_\theta(\bar{\rho},\bar{\theta})}{\grave{o}^2\bar{\theta}}\|\zeta\|_{H^N}^2+\frac{\sigma_a}{4\grave{o}\tilde{\sigma}\bar{\theta}^4}\|\vartheta\|_{H^N}^2\right)
$$

$$
+c_0\|\nabla\mathbf{m}\|_{H^N}^2+\frac{\kappa}{\grave{o}^2\bar{\theta}}\|\nabla\zeta\|_{H^N}^2+\frac{\sigma_a\eta}{4\grave{o}^2\tilde{\sigma}\bar{\theta}^4}\|\nabla\vartheta\|_{H^N}^2+\frac{\sigma_a}{2\grave{o}^2\tilde{\sigma}\bar{\theta}^4}\|(4\tilde{\sigma}\bar{\theta}^3\zeta-\sigma_a\vartheta)\|_{H^N}^2 \quad (2.22)
$$

$$
\le C\|A\|_{H^N}^2\|\mathbf{m}\|_{H^N}^2+C\left(\|\mathbf{G}_1\|_{H^{N-1}}^2+\grave{o}^{-2}\|G_2\|_{H^{N-1}}^2+\grave{o}^{-2}\|G_3\|_{H^{N-1}}^2\right).
$$

On the other hand, noticing that

$$
\int\nabla^k\mathbf{m}_t\cdot\nabla^{k+1}\tilde{n}d\mathbf{x}=\frac{d}{dt}\int\nabla^k\mathbf{m}\cdot\nabla^{k+1}\tilde{n}d\mathbf{x}-\|\nabla^k\mathrm{div}\mathbf{m}\|^2
$$

due to $(2.13)_1$ and multiplying $(2.13)_2$ by $\nabla^{k+1}\tilde{n}$ we have

$$
\frac{d}{dt}\int\nabla^k\mathbf{m}\cdot\nabla^{k+1}\tilde{n}d\mathbf{x}+\frac{P_\rho(\bar{\rho},\bar{\theta})}{\grave{o}^2}\|\nabla^{k+1}\tilde{n}\|^2=J_k^3+\|\nabla^k\mathrm{div}\mathbf{m}\|^2, \quad (2.23)
$$

where

$$
\begin{aligned}
I_k^3=&\frac{P_\theta\left(\bar{\rho},\bar{\theta}\right)}{\grave{o}^2}\int\nabla^k\zeta\cdot\nabla^{k+1}\tilde{n}d\mathbf{x}+\int\nabla^k\mathbf{G_1}\cdot\nabla^{k+1}\tilde{n}d\mathbf{x}\\
&+\bar{\mu}\int\mathrm{div}\nabla^k(A\nabla\mathbf{m})\cdot\nabla^{k+1}\tilde{n}d\mathbf{x}+(\bar{\lambda}+\bar{\mu})\int\nabla\nabla^k(A\mathrm{div}\mathbf{m})\cdot\nabla^{k+1}\tilde{n}d\mathbf{x}.
\end{aligned} \quad (2.24)
$$

Using Holder inequality and (1.18), we obtain from (2.24) that

$$
\begin{aligned}
\sum_{k=0}^{N-1}J_k^3\le& C\left(\grave{o}^{-2}\|\zeta\|_{H^{N-1}}+\|\mathbf{G}_1\|_{H^{N-1}}\right)\|\tilde{n}\|_{H^N}\\
&+C\sum_{k=0}^{N-1}\left(\|A\|_{L^\infty}\|\nabla^{k+2}\mathbf{m}\|+\|\nabla\mathbf{m}\|_{L^\infty}\|\nabla^kA\|\right)\|\nabla^{k+1}\tilde{n}\|\\
\le& C\left(\grave{o}^{-2}\|\zeta\|_{H^N}+\|\mathbf{G}_1\|_{H^{N-1}}\right)\|\tilde{n}\|_{H^N}+C\|A\|_{H^N}\|\nabla\mathbf{m}\|_{H^N}\|\tilde{n}\|_{H^N}.
\end{aligned} \quad (2.25)
$$

Summing up (2.23) for $k=0,1,\cdots,N-1$ and using (2.25), we get

$$
\begin{aligned}
&\frac{d}{dt}\sum_{k=0}^{N-1}\int\nabla^k\mathbf{m}\cdot\nabla^{k+1}\tilde{n}d\mathbf{x}+\frac{P_\rho(\bar{\rho},\bar{\theta})}{\grave{o}^2}\|\nabla\tilde{n}\|_{H^{N-1}}^2\\
&\le\|\nabla\mathbf{m}\|_{H^N}^2+C\left(\grave{o}^{-2}\|\zeta\|_{H^N}+\|\mathbf{G}_1\|_{H^{N-1}}\right)\|\tilde{n}\|_{H^N}+C\|A\|_{H^N}^2\|\tilde{n}\|_{H^N}^2.
\end{aligned} \quad (2.26)
$$

We can assume $0<\grave{o}\le1$ without loss of generality. And we choose $\beta\in(0,1]$ to be suitably small. Then, adding (2.22) and $\beta\times$(2.26), we obtain

$$
\begin{aligned}
&\frac{d}{dt}E(t)+\frac{\beta P_\rho(\bar{\rho},\bar{\theta})}{\grave{o}^2}\|\nabla\tilde{n}\|_{H^{N-1}}^2+(c_0-\beta)\|\nabla\mathbf{m}\|_{H^N}^2\\
&+\frac{\kappa}{\grave{o}^2\bar{\theta}}\|\nabla\zeta\|_{H^N}^2+\frac{\sigma_a}{\eta4\grave{o}^2\tilde{\sigma}\bar{\theta}^4}\|\nabla\vartheta\|_{H^N}^2+\frac{\sigma_a}{2\grave{o}^2\tilde{\sigma}\bar{\theta}^4}\left\|\left(4\tilde{\sigma}\bar{\theta}^3\zeta-\sigma_a\vartheta\right)\right\|_{H^N}^2\\
&\le C\left\|\left(\mathbf{G}_1,\grave{o}^{-1}G_2,\grave{o}^{-1}G_3\right)\right\|_{H^{N-1}}^2+C\left(1+\|A\|_{H^N}^2\right)\left\|\left(\mathbf{m},\grave{o}^{-1}\tilde{n},\grave{o}^{-1}\zeta\right)\right\|_{H^N}^2,
\end{aligned} \quad (2.27)
$$

where

$$E(t) := \|\mathbf{m}\|_{H^N}^2 + \beta \sum_{k=0}^{N-1} \int \nabla^k \mathbf{m} \cdot \nabla^{k+1} \tilde{n} d\mathbf{x} + \frac{\bar{\rho} P_\rho(\bar{\rho}, \bar{\theta})}{\grave{o}^2} \|\tilde{n}\|_{H^N}^2 + \frac{\bar{\rho} e_\theta(\bar{\rho}, \bar{\theta})}{\grave{o}^2 \bar{\theta}} \|\zeta\|_{H^N}^2 + \frac{\sigma_a}{4\grave{o}\tilde{\sigma}\bar{\theta}^4} \| \vartheta \|_{H^N}^2 . \tag{2.28}$$

By (2.28), we can choose the small $\beta \in (0, \frac{c_0}{2}]$ independent $\grave{o} \in (0,1]$ such that

$$E(t) \square \left\| \left( \mathbf{m}, \grave{o}^{-1}\tilde{n}, \grave{o}^{-1}\zeta, \grave{o}^{-\frac{1}{2}}\vartheta \right)(t) \right\|_{H^N}^2 \tag{2.29}$$

uniformly for all $t \in [0,T]$.

Integrating (2.27) over $t \in [0,T]$, and using (2.29), (1.5) and the smalless of $\beta$, we have

$$\begin{aligned} & \left\| \left( \mathbf{m}, \grave{o}^{-1}\tilde{n}, \grave{o}^{-1}\zeta, \grave{o}^{-\frac{1}{2}}\vartheta \right)(t) \right\|_{H^N}^2 + \grave{o}^{-2} \int_0^t \|\nabla \tilde{n}\|_{H^{N-1}}^2 d\tau \\ & + \int_0^t \|\nabla \mathbf{m}\|_{H^N}^2 d\tau + \grave{o}^{-2} \int_0^t \left\| \left( \nabla\zeta, \nabla\vartheta, 4\tilde{\sigma}\bar{\theta}^3\zeta - \sigma_a\vartheta \right) \right\|_{H^N}^2 d\tau \\ & \le \left\| \left( \mathbf{m}, \grave{o}^{-1}\tilde{n}, \grave{o}^{-1}\zeta, \grave{o}^{-\frac{1}{2}}\vartheta \right)(0) \right\|_{H^N}^2 + C\int_0^t \| \left( \mathbf{G}_1, \grave{o}^{-1}G_2, \grave{o}^{-1}G_3 \right) \|_{H^{N-1}}^2 d\tau + \\ & C\int_0^t \left( 1 + \| A \|_{H^N}^2 \right) \| \left( \mathbf{m}, \grave{o}^{-1}\tilde{n}, \grave{o}^{-1}\zeta \right) \|_{H^N}^2 \, d\tau. \end{aligned} \tag{2.30}$$

Applying Gronwall inequality to (2.30) yields

$$\| (\mathbf{m}, \grave{o}^{-1}\tilde{n}, \grave{o}^{-1}\zeta, \grave{o}^{-\frac{1}{2}}\vartheta)(t) \|_{H^N}^2 \le C\Lambda_1(T) e^{C\Lambda_2(T)}, \tag{2.31}$$

where

$$\Lambda_1(T) = \left\| \left( \mathbf{m}, \grave{o}^{-1}\tilde{n}, \grave{o}^{-1}\zeta, \grave{o}^{-\frac{1}{2}}\vartheta \right)(0) \right\|_{H^N}^2 + \int_0^T \left\| \left( \mathbf{G}_1, \grave{o}^{-1}G_2, \grave{o}^{-1}G_3 \right) \right\|_{H^{N-1}}^2 d\tau,$$

$$\Lambda_2(T) = \int_0^T \left( 1 + \|A\|_{H^N}^2 \right) d\tau.$$

By (2.30) and (2.31), we get (2.12). The proof of Theorem 2.1 is completed.

Next, applying Theorem 2.1, it is easy to get the existence and uniqueness of the smooth solution for the linearized system \eqref{324} by the standard methods. We will omit the proof for brevity.

**Theorem 2.2** *Let* $0 < T < \infty$. *Assume that* $(\tilde{n}_0, \mathbf{m}_0, \zeta_0, \vartheta_0) \in H^N(\square^3)$ *for an integer* $N \ge 3$. *Also, suppose that* (2.9) *and* (2.10). *Then, there exists a unique solution* $(\tilde{n}, \mathbf{m}, \zeta, \vartheta)$ *of the linearized system* (2.8) *satisfying* (2.11) *and* (2.12).

## 2.3 Proof of Theorem 1.1

We first prove the local existence of the smooth solutions to the reformulated system (2.3), (2.7) using Theorem 2.1. To this end, we define a set by

$$\mathbf{Z}_M(0,t_0)=\Big\{(\tilde{n},\mathbf{m},\zeta,\vartheta)\in X_N(0,t_0)\Big|\|(\tilde{n},\mathbf{m},\zeta,\vartheta)\|^2_{X_N(0,t_0)}\le M,$$
$$0<m_0^{-1}\le 1+\tilde{n}(\mathbf{x},t),\bar{\theta}+\zeta(\mathbf{x},t)\le m_0\Big\}$$

for a positive constant $m_0>1$ and an integer $N\ge 3$, where

$$M=4C_0\left\|\mathbf{m}_0,\epsilon^{-1}\tilde{n}_0,\epsilon^{-1}\zeta_0,\epsilon^{-\frac{1}{2}}\vartheta_0\right\|^2_{H^N},\ C_0 \text{ is the conxtant determined in (2.12)}, \tag{2.32}$$

$$X_N(0,t_0)=\Big\{(\tilde{n},\mathbf{m},\zeta,\vartheta)\in C\left([0,t_0];H^N(\mathbb{R}^3)\right):\nabla\tilde{n}\in L^2\left(0,t_0;H^{N-1}(\mathbb{R}^3)\right)$$
$$(\nabla\mathbf{m},\nabla\zeta,\nabla\vartheta)\in L^2\left(0,t_0;H^N(\mathbb{R}^3)\right)\Big\},$$

and

$$\begin{aligned}\|(\tilde{n},\mathbf{m},\zeta,\vartheta)\|^2_{X_N(0,t_0)}&=\left\|\left(\mathbf{u},\epsilon^{-1}\tilde{n},\epsilon^{-1}\zeta,\epsilon^{-\frac{1}{2}}\vartheta\right)(t)\right\|^2_{H^N}+\epsilon^{-2}\int_0^{t_0}\|\nabla\tilde{n}\|^2_{H^{N-1}}d\tau\\&\quad+\int_0^{t_0}\|\nabla\mathbf{m}\|^2_{H^N}d\tau+\epsilon^{-2}\int_0^t\|(\nabla\zeta,\nabla\vartheta)\|^2_{H^N}d\tau.\end{aligned}\tag{2.33}$$

For the fixed $(\tilde{n}^{j-1},\mathbf{m}^{j-1},\zeta^{j-1},\vartheta^{j-1})\in\mathbf{Z}_N(0,t_0),\ j=1,2,\cdots,$ we consider the following linearized system

$$\begin{aligned}&\tilde{n}_t+\mathrm{div}\mathbf{m}=0,\qquad \mathbf{x}\in\mathbb{R}^3,t>0,\\&\mathbf{m}_t-\bar{\mu}\mathrm{div}\left(\frac{\nabla\mathbf{m}}{1+\tilde{n}^{j-1}}\right)-(\bar{\lambda}+\bar{\mu})\nabla\left(\frac{\mathrm{div}\mathbf{m}}{1+\tilde{n}^{j-1}}\right)\\&\quad+\frac{\bar{\rho}P_\rho(\bar{\rho},\bar{\theta})}{\epsilon^2}\nabla\tilde{n}+\frac{P_\theta(\bar{\rho},\bar{\theta})}{\epsilon^2}\nabla\zeta=\mathbf{G}_1(\tilde{n}^{j-1},\mathbf{m}^{j-1},\zeta^{j-1}),\\&\zeta_t+\frac{\bar{\theta}P_\theta(\bar{\rho},\bar{\theta})}{\bar{\rho}e_\theta(\bar{\rho},\bar{\theta})}\mathrm{div}\mathbf{m}=\frac{\kappa}{\bar{\rho}e_\theta(\bar{\rho},\bar{\theta})}\Delta\zeta-\frac{4\tilde{\sigma}\bar{\theta}^3\zeta-\sigma_a\vartheta}{\bar{\rho}e_\theta(\bar{\rho},\bar{\theta})}+G_2(\tilde{n}^{j-1},\mathbf{m}^{j-1},\zeta^{j-1}),\\&\epsilon\vartheta_t-\nu\Delta\vartheta=4\tilde{\sigma}\bar{\theta}^3\zeta-\sigma_a\vartheta+G_3(\zeta^{j-1}),\\&(\tilde{n},\mathbf{m},\zeta,\vartheta)(\mathbf{x},0)=(\tilde{n}_0,\mathbf{m}_0,\zeta_0,\vartheta_0)(\mathbf{x}),\quad \mathbf{x}\in\Omega,\end{aligned}\tag{2.34}$$

where $(\tilde{n}^0,\mathbf{m}^0,\zeta^0,\vartheta^0)=(\tilde{n}_0,\mathbf{m}_0,\zeta_0,\vartheta_0)$, and $\mathbf{G}_1(\tilde{n}^{j-1},\mathbf{m}^{j-1},\zeta^{j-1})$, $G_2(\tilde{n}^{j-1},\mathbf{m}^{j-1},\zeta^{j-1})$ and $G_3(\zeta^{j-1})$ are defined by (2.4), (2.5) and (2.6), respectively.

By using Lemma 1.2 and (1.17), we have

$$\begin{aligned}\int_0^{t_0}\left\|\mathrm{div}\left(\frac{\mathbf{m}^{j-1}\otimes\mathbf{m}^{j-1}}{1+\tilde{n}^{j-1}}\right)\right\|^2_{H^{N-1}}d\tau&\le\int_0^{t_0}\left\|\left(\frac{\mathbf{m}^{j-1}\otimes\mathbf{m}^{j-1}}{1+\tilde{n}^{j-1}}\right)\right\|^2_{H^N}d\tau\\&\le C\sup_{0\le t\le t_0}\|\mathbf{m}^{j-1}\|^2_{H^N}(1+\sup_{0\le t\le t_0}\|\tilde{n}^{j-1}\|^2_{H^N}\int_0^{t_0}\|\mathbf{m}^{j-1}\|^2_{H^N}d\tau\\&\le CM^2(1+M)t_0,\end{aligned}$$

$$
\int_0^{t_0}\left\|\operatorname{div}\left(\mathbf{m}^{j-1}\otimes\nabla\left(\frac{1}{1+\tilde{n}^{j-1}}\right)\right)\right\|_{H^{N-1}}^2 d\tau+\int_0^{t_0}\left\|\nabla\left(\mathbf{m}^{j-1}\otimes\nabla\left(\frac{1}{1+\tilde{n}^{j-1}}\right)\right)\right\|_{H^{N-1}}^2 d\tau
$$

$$
\le 2\int_0^{t_0}\left\|\left(\mathbf{m}^{j-1}\otimes\nabla\left(\frac{1}{1+\tilde{n}^{j-1}}\right)\right)\right\|_{H^N}^2 d\tau \le C\sup_{0\le t\le t_0}\left\|\mathbf{m}^{j-1}\right\|_{H^N}^2\sup_{0\le t\le t_0}\left\|\nabla\tilde{n}^{j-1}\right\|_{H^{N-1}}^2 t_0
$$

$$
\le CM^2 t_0
$$

and

$$
\int_0^{t_0}\left\|\frac{h_1(\tilde{n}^{j-1},\zeta^{j-1})}{\grave{o}^2}\nabla\tilde{n}^{j-1}\right\|_{H^{N-1}}^2 d\tau+\int_0^{t_0}\left\|\frac{\bar{\rho}h_2(\tilde{n}^{j-1},\zeta^{j-1})}{\grave{o}^2}\nabla\zeta^{j-1}\right\|_{H^{N-1}}^2 d\tau
$$

$$
\le C\grave{o}^{-4}\sup_{0\le t\le t_0}\|(\tilde{n}^{j-1},\zeta^{j-1})\|_{H^{N-1}}^2\int_0^{t_0}\|(\nabla\tilde{n}^{j-1},\nabla\zeta^{j-1})\|_{H^{N-1}}^2 d\tau\le CM^2 t_0
$$

for $(\tilde{n}^{j-1},\mathbf{m}^{j-1},\zeta^{j-1},\vartheta^{j-1})\in\mathbf{Z}_N(0,t_0),\ j=1,2,\cdots$. Therefore, we obtain from (2.4) that

$$
\begin{aligned}
&\mathbf{G}_1(\tilde{n}^{j-1},\mathbf{m}^{j-1},\zeta^{j-1})\in L^2(0,t_0;H^{N-1}(\square^3)),\\
&\int_0^{t_0}\|\mathbf{G}_1(\tilde{n}^{j-1},\mathbf{m}^{j-1},\zeta^{j-1})\|_{H^{N-1}}^2 d\tau\le C(1+M)^3 t_0.
\end{aligned}
\tag{2.35}
$$

By the similar arguments, we obtain from (2.5) and (2.6), respectively, that

$$
\begin{aligned}
&G_2(\tilde{n}^{j-1},\mathbf{m}^{j-1},\zeta^{j-1})\in L^2(0,t_0;H^{N-1}(\square^3)),\\
&\int_0^{t_0}\|G_2(\tilde{n}^{j-1},\mathbf{m}^{j-1},\zeta^{j-1})\|_{H^{N-1}}^2 d\tau\le C(1+M)^5(1+t_0)
\end{aligned}
\tag{2.36}
$$

and

$$
G_3(\zeta^{j-1})\in L^2(0,t_0;H^{N-1}(\square^3)),\quad \int_0^{t_0}\|G_3(\zeta^{j-1})\|_{H^{N-1}}^2 d\tau\le CM^4 t_0. \tag{2.37}
$$

Also, we have

$$
\frac{1}{1+\tilde{n}^{j-1}}\in L^\infty(0,t_0;H^{N-1}(\square^3)),\quad 0<c_0\le\frac{1}{1+\tilde{n}^{j-1}}\le M_0<\infty. \tag{2.38}
$$

By using (2.35)-(2.38) and Theorem 2.2, there exists a unique solution $(\tilde{n}^j,\mathbf{m}^j,\zeta^j,\vartheta^j)\in X_N(0,t_0)$ to the linearized system (2.34) satisfying

$$
\begin{aligned}
&\left\|(\mathbf{m}^j,\grave{o}^{-1}\tilde{n}^j,\grave{o}^{-1}\zeta^j,\grave{o}^{-1/2}\vartheta^j)(t)\right\|_{H^N}^2+\grave{o}^{-2}\int_0^t\left\|\nabla\tilde{n}^j\right\|_{H^{N-1}}^2 d\tau\\
&+\int_0^t\left\|\nabla\mathbf{m}^j\right\|_{H^N}^2 d\tau+\grave{o}^{-2}\int_0^t\left\|\left(\nabla\zeta^j,\nabla\vartheta^j\right)\right\|_{H^N}^2 d\tau\\
&\le C_0\left[\left\|\mathbf{m}_0,\grave{o}^{-1}\tilde{n}_0,\grave{o}^{-1}\zeta_0,\grave{o}^{-1/2}\vartheta_0\right\|_{H^N}^2+\int_0^t\left\|\left(\mathbf{G}_1^{j-1},\grave{o}^{-1}G_2^{j-1},\grave{o}^{-1}G_3^{j-1}\right)\right\|_{H^{N-1}}^2 d\tau\right]\times\\
&\times\left[1+t_0\left(1+\sup_{0\le t\le t_0}\left\|(1+\tilde{n}^{j-1})^{-1}\right\|_{H^N}^2\right)e^{c_0t_0\left(1+\sup_{0\le t\le t_0}\left\|(1+\tilde{n}^{j-1})^{-1}\right\|_{H^N}^2\right)}\right]
\end{aligned}
\tag{2.39}
$$

for any $t\in[0,t_0]$, where $\mathbf{G}_1^{j-1}=\mathbf{G}_1(\tilde{n}^{j-1},\mathbf{m}^{j-1},\zeta^{j-1})$, $G_2^{j-1}=G_2(\tilde{n}^{j-1},\mathbf{m}^{j-1},\zeta^{j-1})$, $G_3^{j-1}=G_3(\zeta^{j-1})$, and we used

$$\int_0^T \left(1+\left\|(1+\varrho^{j-1})^{-1}\right\|_{H^N}^2\right)d\tau \le t_0\left(1+\sup_{0\le t\le t_0}\left\|(1+\varrho^{j-1})^{-1}\right\|_{H^N}^2\right).$$

Therefore, choosing $t_0$ to be small such that

$$t_0\left(1+\sup_{0\le t\le t_0}\left\|(1+\varrho^{j-1})^{-1}\right\|_{H^N}^2\right)e^{c_0 t_0\left(1+\sup_{0\le t\le t_0}\left\|(1+\varrho^{j-1})^{-1}\right\|_{H^N}^2\right)} \le 1$$

and

$$\int_0^t \left\|\left(\mathbf{G}_1^{j-1},\epsilon^{-1}G_2^{j-1},\epsilon^{-1}G_3^{j-1}\right)\right\|_{H^{N-1}}^2 d\tau \le \left\|\mathbf{m}_0,\epsilon^{-1}\varrho_0,\epsilon^{-1}\zeta_0,\epsilon^{-1/2}\vartheta_0\right\|_{H^N}^2,$$

which is possible due to (2.35)-(2.37), we obtain from (2.39) that

$$\begin{aligned}&\left\|(\mathbf{m}^j,\epsilon^{-1}\varrho^j,\epsilon^{-1}\zeta^j,\epsilon^{-1/2}\vartheta^j)(t)\right\|_{H^N}^2+\epsilon^{-2}\int_0^t\left\|\nabla\varrho^j\right\|_{H^{N-1}}^2 d\tau\\&+\int_0^t\left\|\nabla\mathbf{m}^j\right\|_{H^N}^2 d\tau+\epsilon^{-2}\int_0^t\left\|\left(\nabla\zeta^j,\nabla\vartheta^j\right)\right\|_{H^N}^2 d\tau\le 4C_0\left\|\mathbf{m}_0,\epsilon^{-1}\varrho_0,\epsilon^{-1}\zeta_0,\epsilon^{-1/2}\vartheta_0\right\|_{H^N}^2\end{aligned}\tag{2.40}$$

for any $t\in[0,t_0]$, which implies

$$(\varrho^j,\mathbf{m}^j,\zeta^j,\vartheta^j)\in\mathbf{Z}_N(0,t_0)$$

for $j=1,2,\cdots$ due to (2.32) and (2.33). Moreover, by using (2.40) and (1.8), we have

$$\begin{aligned}&\{(\mathbf{m}^j,\varrho^j,\zeta^j,\vartheta^j)\}_{j=1}^\infty \text{ is bounded in } L^\infty(0,t_0;H^N(\mathbb{R}^3)),\\&\{\nabla\varrho^j\}_{j=1}^\infty \text{ is bounded in } L^2(0,t_0;H^{N-1}(\mathbb{R}^3)),\\&\left\{\left(\nabla\mathbf{m}^j,\nabla\zeta^j,\nabla\vartheta^j\right)\right\}_{j=1}^\infty \text{ is bounded in } L^2(0,t_0;H^N(\mathbb{R}^3))\end{aligned}\tag{2.41}$$

for the fixed $\epsilon>0$. Also, by using (2.41) and (2.35)-(2.38), we obtain from (2.34) that

$$\{\partial_t(\mathbf{m}^j,\varrho^j,\zeta^j,\vartheta^j)\}_{j=1}^\infty \text{ is bounded in } L^2(0,t_0;H^{N-1}(\mathbb{R}^3))\tag{2.42}$$

for the fixed $\epsilon>0$.

By using (2.41), (2.42) and the compactness result, there exists the subsequences of $\{(\mathbf{m}^j,\varrho^j,\zeta^j,\vartheta^j)\}$ (denoting them as $\{(\mathbf{m}^j,\varrho^j,\zeta^j,\vartheta^j)\}$ still) such that when $j\to\infty$, it holds that

$$\begin{aligned}&(\mathbf{m}^j,\varrho^j,\zeta^j,\vartheta^j)\to(\mathbf{m},\varrho,\zeta,\vartheta) \text{ (*-weak) in } L^\infty(0,t_0;H^N(\mathbb{R}^3)),\\&\nabla\varrho^j\to\nabla\varrho \text{ (weak) in } L^2(0,t_0;H^{N-1}(\mathbb{R}^3)),\\&\left(\nabla\mathbf{m}^j,\nabla\zeta^j,\nabla\vartheta^j\right)\to\left(\nabla\mathbf{m},\nabla\zeta,\nabla\vartheta\right) \text{ (weak) in } L^2(0,t_0;H^N(\mathbb{R}^3)),\\&(\mathbf{m}^j,\varrho^j,\zeta^j,\vartheta^j)\to(\mathbf{m},\varrho,\zeta,\vartheta) \text{ (strong) in } L^2(0,t_0;H_{loc}^{N-1}(\mathbb{R}^3)),\end{aligned}\tag{2.43}$$

where $H_{loc}^{N-1}(\mathbb{R}^3)$ denote the space $\{w\ |w\in H^{N-1}(\Omega)\}$ for any bounded domain $\Omega$ of $\mathbb{R}^3$, and

$$\begin{aligned}&(\mathbf{m},\varrho,\zeta,\vartheta)\in L^\infty(0,t_0;H^N(\mathbb{R}^3)),\\&\nabla\varrho\in L^2(0,t_0;H^{N-1}(\mathbb{R}^3)),\quad\left(\nabla\mathbf{m},\nabla\zeta,\nabla\vartheta\right)\in L^2(0,t_0;H^N(\mathbb{R}^3)).\end{aligned}\tag{2.44}$$

By using (2.43) and (2.4)-(2.6), we have

$$\int_0^{t_0} \mathbf{G}_1(\tilde{n}^{j-1}, \mathbf{m}^{j-1}, \zeta^{j-1}) \cdot \mathbf{w}\,\eta(t) d\mathbf{x}dt \to \int_0^{t_0} \mathbf{G}_1(\tilde{n}, \mathbf{m}, \zeta) \cdot \mathbf{w}\,\eta(t) d\mathbf{x}dt,$$
$$\int_0^{t_0} G_2(\tilde{n}^{j-1}, \mathbf{m}^{j-1}, \zeta^{j-1}) w\eta(t) d\mathbf{x}dt \to \int_0^{t_0} \mathbf{G}_1(\tilde{n}, \mathbf{m}, \zeta) \cdot \mathbf{w}\,\eta(t) d\mathbf{x}dt, \qquad (2.45)$$
$$\int_0^{t_0} G_3(\zeta^{j-1}) w\eta(t) d\mathbf{x}dt \to \int_0^{t_0} G_3(\zeta) w\eta(t) d\mathbf{x}dt$$

for all $\mathbf{w} \in C_0^\infty(\mathbb{R}^3)^3$, $w \in C_0^\infty(\mathbb{R}^3)$ and $\eta \in C_0^\infty(0, t_0)$. Moreover, by using (2.43) and (2.45), it is easy to check that $(\mathbf{m}, \tilde{n}, \zeta, \vartheta)$ is a solution to the system (2.3), (2.7). Then, setting

$$\rho = \bar{\rho}(1+\tilde{n}), \quad \mathbf{u} = \frac{\mathbf{m}}{1+\tilde{n}}, \quad \theta = \bar{\theta} + \zeta \quad \text{and} \quad n = \bar{n} + \vartheta,$$

and by using (2.3), (2.7), (2.4)-(2.6) and (2.44), we know that $(\rho, \mathbf{u}, \theta, n)$ is a solution to the Cauchy problem (1.3), (1.4) satisfying (1.9). Also, the estimate (1.10) follows from (2.40) and (2.43). The proof for the uniqueness of the local solution is standard, so we will omit it for brevity.
We prove the low Mach number limit. Let the solution to the Cauchy problem (1.3), (1.4) satisfying (1.9) and (1.10) by $(\rho^\varepsilon, \mathbf{u}^\varepsilon, \theta^\varepsilon, n^\varepsilon)$. Then, using (1.10), Lemma 1.2, $(1.3)_{1,2}$ and (2.1), there exists a function $\mathbf{u} \in L^\infty(0, T_0; H^N(\mathbb{R}^3)) \cap L^2(0, T_0; H^{N+1}(\mathbb{R}^3))$ such that

$$\begin{aligned}
&(\rho^\varepsilon, \mathbf{u}^\varepsilon, \theta^\varepsilon, n^\varepsilon) \to (\bar{\rho}, \mathbf{u}, \bar{\theta}, \bar{n}) \text{ (*-weak) in } L^\infty(0, t_0; H^N(\mathbb{R}^3)),\\
&\nabla \tilde{n}^\varepsilon \to 0 \text{ (weak) in } L^2(0, t_0; H^{N-1}(\mathbb{R}^3)),\\
&\left(\nabla \mathbf{u}^\varepsilon, \nabla \zeta^\varepsilon, \nabla \vartheta^\varepsilon\right) \to \left(\nabla \mathbf{u}, 0, 0\right) \text{ (weak) in } L^2(0, t_0; H^N(\mathbb{R}^3)),\\
&(\rho^\varepsilon, \mathbf{u}^\varepsilon, \theta^\varepsilon) \to (\bar{\rho}, \mathbf{u}, \bar{\theta}) \text{ (strong) in } L^2(0, t_0; H_{loc}^{N-1}(\mathbb{R}^3)),
\end{aligned} \qquad (2.46)$$

as $\varepsilon \to 0$, which implies (1.11).
Now, we will prove that the function pair $(\mathbf{u}, P)$ for some $P \in L^\infty(0, T_0; H^{N-1}(\mathbb{R}^3))$ $\bigcap L^2(0, T_0; H^N(\mathbb{R}^3))$ is the unique smooth local solution for the incompressible Navier-Stokes system (1.12), (1.13). By $(1.3)_1$ and (2.46), we get $\mathrm{div}\mathbf{u} = 0$, which implies $(1.12)_1$. Multiplying $(1.3)_2$ by $\mathbf{w} \in [C_0^\infty(\mathbb{R}^3)]^3, \mathrm{div}\mathbf{w} = 0,$ and integrating it over $\mathbb{R}^3$, we have

$$\int \mathbf{u}_t^\varepsilon \cdot \mathbf{w} d\mathbf{x} + \int (\mathbf{u}^\varepsilon \cdot \nabla)\mathbf{u}^\varepsilon \cdot \mathbf{w} d\mathbf{x} = \int \mu \rho^\varepsilon \Delta \mathbf{u}^\varepsilon \cdot \mathbf{w} d\mathbf{x}. \qquad (2.47)$$

By using (2.46) and the result of De Rham ([18, Proposition 1.1]), we obtain $(1.12)_2$ from (2.47). It is easy to show that the initial condition (1.13) holds. The uniqueness of solutions to the system (1.12), (1.13) is proved by a standard argument.
The proof of Theorem 1.1 is completed.

## 3. Global solution and low Mach number limit

In this section, we prove Theorem 1.2. For brevity, we denote $(\rho^\varepsilon, \mathbf{u}^\varepsilon, \theta^\varepsilon, n^\varepsilon)$, $P^\varepsilon = P(\rho^\varepsilon, \theta^\varepsilon)$ and $e^\varepsilon = e(\rho^\varepsilon, \theta^\varepsilon)$ in (1.3 ) by $(\rho, \mathbf{u}, \theta, n)$, $P = P(\rho, \theta)$ and $e = e(\rho, \theta)$, respectively.

### 3.1 A priori estimates

This subsection is devoted to derive the a priori estimates, which play a crucial role in the proof of Theorem 1.2. Setting

$$\varphi = \rho - \bar{\rho}, \quad \zeta = \theta - \bar{\theta} \quad \text{and} \quad \vartheta = n - \bar{n},$$

Therefore, using (2.1) and (2.2), we reformulate the system (1.3) as

$$\begin{aligned}
&\varphi_t + \bar{\rho}\mathrm{div}\mathbf{u} = g_1,\\
&\mathbf{u}_t - \frac{\mu}{\bar{\rho}}\Delta\mathbf{u} - \frac{\mu+\lambda}{\bar{\rho}}\nabla\mathrm{div}\mathbf{u} + \frac{P_\rho(\bar{\rho},\bar{\theta})}{\bar{\rho}\epsilon^2}\nabla\varphi + \frac{P_\theta(\bar{\rho},\bar{\theta})}{\bar{\rho}\epsilon^2}\nabla\zeta = \mathbf{g}_2,\\
&\zeta_t + \frac{\bar{\theta}P_\theta(\bar{\rho},\bar{\theta})}{e_\theta(\bar{\rho},\bar{\theta})}\mathrm{div}\mathbf{u} = \frac{\kappa}{\bar{\rho}e_\theta(\bar{\rho},\bar{\theta})}\Delta\theta - \frac{4\tilde{\sigma}\bar{\theta}^3\zeta - \sigma_a\vartheta}{\bar{\rho}e_\theta(\bar{\rho},\bar{\theta})} + g_3,\\
&\epsilon n_t - \nu\Delta n = 4\tilde{\sigma}\bar{\theta}^3\zeta - \sigma_a\vartheta + g_4,
\end{aligned} \tag{3.1}$$

where $g_1, \mathbf{g}_2, g_3$ and $g_4$ are defined respectively by

$$g_1 = -\varphi\mathrm{div}\mathbf{u} - \mathbf{u}\cdot\nabla\varphi, \tag{3.2}$$

$$\begin{aligned}
&\mathbf{g}_2 = -(\mathbf{u},\nabla)\mathbf{u} + \frac{h_6(\varphi,\zeta)}{\epsilon^2}\nabla\varphi + \frac{h_7(\varphi,\zeta)}{\epsilon^2}\nabla\zeta - h_8(\varphi)(\mu\Delta\mathbf{u} + (\mu+\lambda)\nabla\mathrm{div}\mathbf{u}),\\
&h_6(\varphi,\zeta) = \frac{P_\rho(\bar{\rho},\bar{\theta})}{\bar{\rho}} - \frac{P_\rho(\rho,\theta)}{\rho}, \; h_7(\varphi,\zeta) = \frac{P_\theta(\bar{\rho},\bar{\theta})}{\bar{\rho}} - \frac{P_\theta(\rho,\theta)}{\rho},\\
&h_8(\varphi) = \bar{\rho}^{-1} - \rho^{-1},
\end{aligned} \tag{3.3}$$

$$\begin{aligned}
g_3 = &-\mathbf{u}\cdot\nabla\zeta - \kappa h_9(\varphi,\zeta)\Delta\zeta + \frac{2\mu\epsilon^2\mathbf{D}(\mathbf{u}):\mathbf{D}(\mathbf{u})}{\rho e_\theta(\rho,\theta)} + \frac{\lambda\epsilon^2(\mathrm{div}\mathbf{u})^2}{\rho e_\theta(\rho,\theta)}\\
&+ h_{10}(\varphi,\zeta)\mathrm{div}\mathbf{u} - h_9(\varphi,\zeta)(4\tilde{\sigma}\bar{\theta}^3\zeta - \sigma_a\vartheta) - \frac{h_5(\zeta)\zeta}{\rho e_\theta(\rho,\theta)},\\
&h_9(\varphi,\zeta) = \frac{1}{\bar{\rho}e_\theta(\bar{\rho},\bar{\theta})} - \frac{1}{\rho e_\theta(\rho,\theta)}, \; h_{10}(\varphi,\zeta) = \frac{\bar{\theta}P_\theta(\bar{\rho},\bar{\theta})}{\bar{\rho}e_\theta(\bar{\rho},\bar{\theta})} - \frac{\theta P_\theta(\rho,\theta)}{\rho e_\theta(\rho,\theta)},\\
&h_5(\zeta) = 6\tilde{\sigma}\bar{\theta}^2\zeta + 4\tilde{\sigma}\bar{\theta}\zeta^2 + \tilde{\sigma}\zeta^3
\end{aligned} \tag{3.4}$$

and

$$g_4 = h_5(\zeta)\zeta. \tag{3.5}$$

We suppose that $(\varphi,\mathbf{u},\zeta,\vartheta) \in X_N(0,T)$ is the solution to the system (3.1) for any constant $T > 0$ and an integer $N \geq 3$, where $X_N(I)$ is defined by

$$\begin{aligned}
X_N(I) = \{&(\varphi,\mathbf{u},\zeta,\vartheta) \in C(I;H^N(\mathbb{R}^3)) \mid \nabla\varphi \in L^2(I;H^{N-1}(\mathbb{R}^3)),\\
&(\nabla\mathbf{u},\nabla\zeta,\nabla\vartheta,4\tilde{\sigma}\bar{\theta}^3\zeta - \sigma_a\vartheta) \in L^2(I;H^N(\mathbb{R}^3))\}
\end{aligned}$$

for any interval $I \subset [0,\infty)$. Also, we assume

$$\sup_{0\leq t\leq T}\left(\|\mathbf{u}(t)\|_{H^3} + \epsilon^{-1}\|(\varphi,\zeta)(t)\|_{H^3} + \epsilon^{-\frac{1}{2}}\|\vartheta(t)\|_{H^3}\right) \leq \delta_1 \tag{3.6}$$

for sufficiently small $\delta_1 > 0$ independent of $\epsilon$. By using (3.6), we can choose $\delta_1 > 0$ such that

$$0 < \frac{\bar{\rho}}{2} \leq \rho(x,t) \leq 2\bar{\rho} \quad \text{and} \quad 0 < \frac{\bar{\theta}}{2} \leq \theta(x,t) \leq 2\bar{\theta}. \tag{3.7}$$

We first give the following estimate:

**Lemma 3.1** *Under the assumption* (3.6), *it holds that*

$$\frac{d}{dt}\left(\|\mathbf{u}\|_{H^l}^2+\frac{P_\rho(\bar\rho,\bar\theta)}{\bar\rho^2ò^2}\|\varphi\|_{H^l}^2+\frac{e_\theta(\bar\rho,\bar\theta)}{\bar\theta ò^2}\|\zeta\|_{H^l}^2+\frac{\sigma_a}{4\tilde\sigma ò\bar\rho\bar\theta^4}\|\vartheta\|_{H^l}^2\right)$$
$$+\frac{\mu}{\bar\rho}\|\nabla\mathbf{u}\|_{H^l}^2+\frac{\kappa}{\bar\rho\bar\theta ò^2}\|\nabla\zeta\|_{H^l}^2+\frac{\nu\sigma_a}{4\tilde\sigma\bar\rho\bar\theta^4ò^2}\|\nabla\vartheta\|_{H^l}^2 \tag{3.8}$$
$$+\frac{1}{4\tilde\sigma\bar\rho\bar\theta^4ò^2}\|(4\tilde\sigma\bar\theta^3\zeta-\sigma_a\vartheta)\|_{H^l}^2\le\frac{C\delta_1}{ò^2}\|\nabla\varphi\|_{H^{l-1}}^2$$

*for any* $t\ge 0$ *and* $l=3,\cdots,N$, *where* $C$ *is the positive constant independent of* $ò,\delta_1$ *and* $t$.

**Proof.** Applying $\nabla^k$ to (3.1) yields

$$\nabla^k\varphi_t+\bar\rho\,\mathrm{div}\nabla^k\mathbf{u}=\nabla^k g_1,$$
$$\nabla^k\mathbf{u}_t-\frac{\mu}{\bar\rho}\Delta\nabla^k\mathbf{u}-\frac{(\mu+\lambda)}{\bar\rho}\nabla^{k+1}\mathrm{div}\mathbf{u}+\frac{P_\rho(\bar\rho,\bar\theta)}{\bar\rho ò^2}\nabla^{k+1}\varphi+\frac{P_\theta(\bar\rho,\bar\theta)ò^2}{\bar\rho}\nabla^{k+1}\zeta=\nabla^k\mathbf{g}_2, \tag{3.9}$$
$$\nabla^k\zeta_t+\frac{\bar\theta P_\theta(\bar\rho,\bar\theta)}{e_\theta(\bar\rho,\bar\theta)}\mathrm{div}\nabla^k\mathbf{u}=\frac{\kappa}{\bar\rho e_\theta(\bar\rho,\bar\theta)}\Delta\nabla^k\theta+\frac{\nabla^k(4\tilde\sigma\bar\theta^3\zeta-\sigma_a\vartheta)}{\bar\rho e_\theta(\bar\rho,\bar\theta)}+\nabla^k g_3,ò\nabla$$
$${}^k\vartheta_t-\nu\nabla^k\Delta\vartheta=\nabla^k(4\tilde\sigma\bar\theta^3\zeta-\sigma_a\vartheta)+\nabla^k g_4,$$

where $k=0,1,\cdots,l.$

Multiplying (3.9)2, (3.9)3 and (3.9)4 by $\nabla^k\mathbf{u},\frac{e_\theta(\bar\rho,\bar\theta)}{ò^2\bar\theta}\nabla^k\zeta$ and $\frac{\sigma_a}{4ò^2\tilde\sigma\bar\rho\bar\theta^4}\nabla^k\vartheta$, respectively, and using (3.9)1 and (3.2)-(35), we have

$$\frac{1}{2}\frac{d}{dt}\left(\|\nabla^k\mathbf{u}\|^2+\frac{P_\rho(\bar\rho,\bar\theta)}{\bar\rho^2ò^2}\|\nabla^k\varphi\|^2+\frac{e_\theta(\bar\rho,\bar\theta ò^2)}{\bar\theta}\|\nabla^k\zeta\|^2+\frac{\sigma_a}{4ò\tilde\sigma\bar\rho\bar\theta^4}\|\nabla^k\vartheta\|^2\right)$$
$$+\int\left(\frac{\mu}{\bar\rho}|\nabla^{k+1}\mathbf{u}|^2+\frac{(\mu+\lambda)}{\bar\rho}|\mathrm{div}\nabla^k\mathbf{u}|^2\right)d\mathbf{x}+\frac{\kappa}{\bar\rho\bar\theta ò^2}\int|\nabla^{k+1}\zeta|^2\,d\mathbf{x} \tag{3.10}$$
$$+\frac{\nu\sigma_a}{4ò^2\tilde\sigma\bar\rho\bar\theta^4}\int|\nabla^{k+1}\vartheta|^2\,d\mathbf{x}+\frac{1}{4ò^2\tilde\sigma\bar\rho\bar\theta^4}\int|\nabla^k(4\tilde\sigma\bar\theta^3\zeta-\sigma_a\vartheta)|^2\,d\mathbf{x}$$
$$=\frac{P_\theta(\bar\rho,\bar\theta)}{\bar\rho^2}I_k^1+I_k^2+\frac{e_\theta(\bar\rho,\bar\theta)}{\bar\theta}I_k^3+\frac{1}{4\tilde\sigma\bar\rho\bar\theta^4}I_k^4,$$

where

$$I_k^1=-ò^{-2}\int\nabla^k(\varphi\mathrm{div}\mathbf{u}+\mathbf{u}\nabla\varphi)\nabla^k\varphi d\mathbf{x}, \tag{3.11}$$

$$I_k^2=\underbrace{-\int\nabla^k[(\mathbf{u},\nabla)\mathbf{u}]\cdot\nabla^k\mathbf{u}d\mathbf{x}}_{I_k^{21}}+\underbrace{ò^{-2}\int\nabla^k[h_6(\varphi,\zeta)\nabla\varphi+h_7(\varphi,\zeta)\nabla\zeta]\cdot\nabla^k\mathbf{u}d\mathbf{x}}_{I_k^{22}}$$
$$\underbrace{-\int\nabla^k[h_8(\varphi)(\mu\Delta\mathbf{u}+(\mu+\lambda)\nabla\mathrm{div}\mathbf{u})]\cdot\nabla^k\mathbf{u}d\mathbf{x}}_{I_k^{23}}, \tag{3.12}$$

$$
\begin{aligned}
I_k^3 = & \underbrace{-\kappa\epsilon^{-2}\int \nabla^k[h_9(\varphi,\zeta)\Delta\zeta]\nabla^k\zeta d\mathbf{x}}_{I_k^{31}} \\
& \underbrace{-\epsilon^{-2}\int \nabla^k[\mathbf{u}\cdot\nabla\zeta - h_{10}(\varphi,\zeta)\mathrm{div}\mathbf{u}]\nabla^k\zeta d\mathbf{x}}_{I_k^{32}} \\
& \underbrace{+\int \nabla^k\left[\frac{2\mu\mathbf{D}(\mathbf{u}):\mathbf{D}(\mathbf{u})+\lambda(\mathrm{div}\mathbf{u})^2}{\rho e_\theta(\rho,\theta)}\right]\nabla^k\zeta d\mathbf{x}}_{I_k^{33}} \\
& \underbrace{-\epsilon^{-2}\int \nabla^k[h_9(\varphi,\zeta)(4\tilde{\sigma}\bar{\theta}^3\zeta-\sigma_a\vartheta)]\nabla^k\zeta d\mathbf{x}}_{I_k^{34}} \underbrace{-\epsilon^{-2}\int \nabla^k[h_9(\varphi,\zeta)h_5(\zeta)\zeta]\nabla^k\zeta d\mathbf{x}}_{I_k^{35}}
\end{aligned}
\tag{3.13}
$$

and

$$
I_k^4 = -\epsilon^{-2}\int \nabla^k[h_5(\zeta)\zeta]\nabla^k(4\tilde{\sigma}\bar{\theta}^3\zeta-\sigma_a\vartheta)d\mathbf{x}. \tag{3.14}
$$

Adding (3.10) for $k=0,1,\cdots,l$ yields that

$$
\begin{aligned}
& \frac{1}{2}\frac{d}{dt}\Big(\|\mathbf{u}\|_{H^l}^2 + \frac{P_\rho(\bar{\rho},\bar{\theta})}{\bar{\rho}^2\epsilon^2}\|\varphi\|_{H^l}^2 + \frac{e_\theta(\bar{\rho},\bar{\theta})}{\bar{\theta}\epsilon^2}\|\zeta\|_{H^l}^2 + \frac{\sigma_a}{4\tilde{\sigma}\epsilon\bar{\rho}\bar{\theta}^4}\|\vartheta\|_{H^l}^2\Big) \\
& + \frac{\mu}{\bar{\rho}}\|\nabla\mathbf{u}\|_{H^l}^2 + \frac{\kappa}{\bar{\rho}\bar{\theta}\epsilon^2}\|\nabla\zeta\|_{H^l}^2 + \frac{\nu\sigma_a}{4\tilde{\sigma}\bar{\rho}\bar{\theta}^4\epsilon^2}\|\nabla\vartheta\|_{H^l}^2 + \frac{1}{4\tilde{\sigma}\bar{\rho}\bar{\theta}^4\epsilon^2}\|(4\tilde{\sigma}\bar{\theta}^3\zeta-\sigma_a\vartheta)\|_{H^l}^2 \\
& \le \frac{P_\theta(\bar{\rho},\bar{\theta})}{\bar{\rho}^2}\sum_{k=0}^{l}|I_k^1| + \sum_{k=0}^{l}|I_k^2| + \frac{e_\theta(\bar{\rho},\bar{\theta})}{\bar{\theta}}\sum_{k=0}^{l}|I_k^3| + \frac{1}{4\tilde{\sigma}\bar{\rho}\bar{\theta}^4}\sum_{k=0}^{l}|I_k^4|.
\end{aligned}
\tag{3.15}
$$

We estimate the right hands in (3.15).

For $I_0^1$, using Holder inequality, (3.7) and (1.16), we obtain from (3.11) that

$$
I_0^1 \le C\epsilon^{-2}(\|\varphi\|_{L^6}\|\nabla\mathbf{u}\| + \|\mathbf{u}\|_{L^6}\|\nabla\varphi\|)\|\varphi\|_{L^3} \le C\delta_1\epsilon^{-2}\|\nabla\varphi\|^2.
$$

For $I_k^1(k=1,\cdots,l)$, using Holder inequality, (3.7), (1.17)$_2$ and (3.6), we obtain from (3.11) that

$$
\begin{aligned}
I_k^1 &= \epsilon^{-2}\int \nabla^k(\varphi\mathrm{div}\mathbf{u})\nabla^k\varphi d\mathbf{x} - \epsilon^{-2}\int [\nabla^k,\mathbf{u}]\nabla\varphi\nabla^k\varphi d\mathbf{x} + \frac{1}{2\epsilon^2}\int \mathrm{div}\mathbf{u}|\nabla^k\varphi|^2 d\mathbf{x} \\
&\le C\epsilon^{-2}\Big(\|\varphi\|_{L^\infty}\|\nabla^{k+1}\mathbf{u}\| + \|\nabla\mathbf{u}\|_{L^\infty}\|\nabla^k\varphi\|\Big)\|\nabla^k\varphi\| \\
&\quad + C\epsilon^{-2}\Big(\|\nabla\mathbf{u}\|_{L^\infty}\|\nabla^k\varphi\| + \|\nabla\varphi\|_{L^3}\|\nabla^k\mathbf{u}\|_{L^6}\Big)\|\nabla^k\varphi\| + C\epsilon^{-2}\|\nabla\mathbf{u}\|_{L^\infty}\|\nabla^k\varphi\|^2 \\
&\le C\delta_1(\epsilon^{-2}\|\nabla^k\varphi\|^2 + \|\nabla^{k+1}\mathbf{u}\|^2).
\end{aligned}
$$

Therefore, we have

$$
\sum_{k=0}^{l}|I_k^1| \le C\delta_1(\epsilon^{-2}\|\nabla\varphi\|_{H^{l-1}}^2 + \|\nabla\mathbf{u}\|_{H^l}^2). \tag{3.16}
$$

For $I_0^2$, using Holder inequality, (1.16), (3.7) and (3.6), we obtain from (3.12) that

$$I_0^{21} \le C\|\mathbf{u}\|_{L^6}\|\nabla\mathbf{u}\|\|\mathbf{u}\|_{L^3} \le C\delta_1\|\nabla\mathbf{u}\|^2,$$

$$I_0^{22} \le C\epsilon^{-2}(\|h_6(\varphi,\zeta)\|_{L^3}\|\nabla\varphi\| + \|h_7(\varphi,\zeta)\|_{L^3}\|\nabla\zeta\|)\|\mathbf{u}\|_{L^6} \le$$

$$C\delta_1(\|\nabla\mathbf{u}\|^2 + \epsilon^{-2}\|\nabla\varphi\|^2 + \epsilon^{-2}\|\nabla\zeta\|^2),$$

$$I_0^{23} = -\int[\mu\nabla\mathbf{u}\cdot\nabla(h_8(\varphi)\mathbf{u}) + (\lambda+\mu)\mathrm{div}\mathbf{u}\,\mathrm{div}(h_8(\varphi)\mathbf{u})]d\mathbf{x}$$

$$\le C\|\nabla\mathbf{u}\|(\|h_8(\varphi)\|_{L^\infty}\|\nabla\mathbf{u}\| + \|\mathbf{u}\|_{L^\infty}\|h_8(\varphi)\|) \le C\delta_1(\|\nabla\mathbf{u}\|^2 + \|\nabla\varphi\|^2).$$

Therefore, we have

$$I_0^2 \le C\delta_1(\|\nabla\mathbf{u}\|^2 + \epsilon^{-2}\|\nabla\varphi\|^2 + \epsilon^{-2}\|\nabla\zeta\|^2), \tag{3.17}$$

where we assume $0<\epsilon\le 1$ without loss of generality. For $I_k^2 (k=1,\cdots,l)$, we obtain from (3.12) that

$$I_k^{21} \overset{(1.18)}{\le} C\Big(\|\mathbf{u}\|_{L^3}\|\nabla^{k+1}\mathbf{u}\| + \|\nabla\mathbf{u}\|_{L^3}\|\nabla^k\mathbf{u}\|\Big)\|\nabla^k\mathbf{u}\|_{L^6}$$

$$\overset{(3.6)}{\le} C\delta_1\left(\|\nabla^{k+1}\mathbf{u}\|^2 + \|\nabla^k\mathbf{u}\|^2\right),$$

$$I_k^{22} = -\epsilon^{-2}\int\nabla^{k-1}[h_6(\varphi,\zeta)\nabla\varphi + h_7(\varphi,\zeta)\nabla\zeta]\cdot\nabla^{k+1}\mathbf{u}d\mathbf{x}$$

$$\overset{(1.18)}{\le} C\epsilon^{-2}\Big(\|h_6(\varphi,\zeta)\|_{L^\infty}\|\nabla^k\varphi\| + \|\nabla\varphi\|_{L^3}\|\nabla^{k-1}h_6(\varphi,\zeta)\|_{L^6}\Big)\|\nabla^{k+1}\mathbf{u}\|$$

$$+C\epsilon^{-2}\Big(\|h_7(\varphi,\zeta)\|_{L^\infty}\|\nabla^k\zeta\| + \|\nabla\zeta\|_{L^{\frac{3}{2}}}\|\nabla^{k-1}h_7(\varphi,\zeta)\|_{L^6}\Big)\|\nabla^{k+1}\mathbf{u}\|$$

$$\overset{(3.6),(1.17)}{\le} C\delta_1\left(\|\nabla^{k+1}\mathbf{u}\|^2 + \epsilon^{-2}\|\nabla^k\varphi\|^2 + \epsilon^{-2}\|\nabla^k\zeta\|^2\right),$$

$$I_k^{23} \overset{(1.18)}{\le} C\Big(\|h_8(\varphi)\|_{L^\infty}\|\nabla^{k+1}\mathbf{u}\| + \|\nabla\mathbf{u}\|_{L^\infty}\|\nabla^k h_8(\varphi)\| + \|\nabla h_8(\varphi)\|_{L^3}\|\nabla^k\mathbf{u}\|_{L^6}\Big)\|\nabla^{k+1}\mathbf{u}\|$$

$$\overset{(3.6),(1.17)}{\le} C\delta_1\left(\|\nabla^{k+1}\mathbf{u}\|^2 + \|\nabla^k\varphi\|^2\right).$$

Therefore, we have

$$I_k^2 \le C\delta_1\left(\|\nabla^{k+1}\mathbf{u}\|^2 + \epsilon^{-2}\|\nabla^k\varphi\|^2 + \epsilon^{-2}\|\nabla^k\zeta\|^2\right), \tag{3.18}$$

assuming $0<\epsilon\le 1$ without loss of generality. By (3.17) and (3.18), we get

$$\sum_{k=0}^{l}|I_k^2| \le C\delta_1\left(\|\nabla\mathbf{u}\|_{H^l}^2 + \epsilon^{-2}\|\nabla\varphi\|_{H^{l-1}}^2 + \epsilon^{-2}\|\nabla\zeta\|_{H^{l-1}}^2\right). \tag{3.19}$$

For $I_0^3$, using Holder inequality, (3.6) and (1.16), we obtain from (3.13) that

$$I_0^{31} = \epsilon^{-2}\kappa\int\nabla\zeta\cdot\nabla(h_9(\varphi,\zeta)\zeta)d\mathbf{x}$$

$$\le C\epsilon^{-2}\|\nabla\zeta\|\Big(\|h_9(\varphi,\zeta)\|_{L^\infty}\|\nabla\zeta\| + \|\zeta\|_{L^\infty}\|\nabla h_9(\varphi,\zeta)\|\Big)$$

$$\le C\delta_1\epsilon^{-2}\left(\|\nabla\zeta\|^2 + \|\nabla\varphi\|^2\right),$$

$$
\begin{aligned}
I_0^{32} &\le C\epsilon^{-2}\left(\|\mathbf{u}\|_{L^3}\|\nabla\zeta\|+\|h_{10}(\varphi,\zeta)\|_{L^3}\|\nabla\mathbf{u}\|\right)\|\zeta\|_{L^6}\\
&\le C\delta_1\left(\epsilon^{-2}\|\nabla\zeta\|^2+\|\nabla\mathbf{u}\|^2\right),\\
I_0^{33} &\le C\left(\|\mathbf{D}(\mathbf{u})\|_{L^3}\|\mathbf{D}(\mathbf{u})\|+\|\mathrm{div}\mathbf{u}\|_{L^3}\|\mathrm{div}\mathbf{u}\|\right)\|\zeta\|_{L^6}\\
&\le C\delta_1\left(\|\nabla\mathbf{u}\|^2+\|\nabla\zeta\|^2\right),\\
I_0^{34} &\le C\epsilon^{-2}\|h_9(\varphi,\zeta)\|_{L^3}\left\|(4\tilde{\sigma}\bar{\theta}^3\zeta-\sigma_a\vartheta)\right\|\|\zeta\|_{L^6}\\
&\le C\delta_1\epsilon^{-2}\left(\|\nabla\zeta\|^2+\left\|(4\tilde{\sigma}\bar{\theta}^3\zeta-\sigma_a\vartheta)\right\|^2\right),\\
I_0^{35} &\le C\epsilon^{-2}\|h_9(\varphi,\zeta)\|_{L^3}\|h_6(\zeta)\|_{L^3}\|\zeta\|_{L^6}^2\le C\delta_1\epsilon^{-2}\|\nabla\zeta\|^2.
\end{aligned}
$$

Therefore, we have

$$
I_0^3\le C\delta_1\left(\|\nabla\mathbf{u}\|^2+\epsilon^{-2}\|\nabla\varphi\|^2+\epsilon^{-2}\|\nabla\zeta\|^2+\epsilon^{-2}\left\|(4\tilde{\sigma}\bar{\theta}^3\zeta-\sigma_a\vartheta)\right\|^2\right). \tag{3.20}
$$

For $I_k^3(k=1,\cdots,l)$, we obtain from (3.13) that

$$
\begin{aligned}
I_k^{31} &= \kappa\epsilon^{-2}\int\nabla^{k-1}\left[\mathrm{div}(h_9(\varphi,\zeta)\nabla\zeta)-\nabla h_9(\varphi,\zeta)\nabla\zeta\right]\cdot\nabla^{k+1}\zeta\,d\mathbf{x}\\
&\overset{(1.18)}{\le} C\epsilon^{-2}\Big(\|h_9(\varphi,\zeta)\|_{L^\infty}\|\nabla^{k+1}\zeta\|+\|\nabla\zeta\|_{L^\infty}\|\nabla^k h_9(\varphi,\zeta)\|\\
&\qquad+\|\nabla h_4(\varphi,\zeta)\|_{L^3}\|\nabla^k\zeta\|_{L^6}\Big)\|\nabla^{k+1}\zeta\|\\
&\overset{(3.6),(1.17)}{\le} C\delta_1\epsilon^{-2}\left(\|\nabla^{k+1}\zeta\|^2+\|\nabla^k\zeta\|^2+\|\nabla^k\varphi\|^2\right),\\
I_k^{32} &\overset{(1.18)}{\le} C\epsilon^{-2}\Big(\|\mathbf{u}\|_{L^3}\|\nabla^{k+1}\zeta\|+\|\nabla\zeta\|_{L^3}\|\nabla^k\mathbf{u}\|\Big)\|\nabla^k\zeta\|_{L^6}\\
&\quad+\Big(\|h_{10}(\varphi,\zeta)\|_{L^3}\|\nabla^{k+1}\mathbf{u}\|+\|\nabla\mathbf{u}\|_{L^3}\|\nabla^k h_{10}(\varphi,\zeta)\|\Big)\|\nabla^k\zeta\|_{L^6}\\
&\overset{(3.6),(1.17)}{\le} C\delta_1\left(\|\nabla^{k+1}\mathbf{u}\|^2+\|\nabla^k\mathbf{u}\|^2+\epsilon^{-2}\|\nabla^{k+1}\zeta\|^2+\epsilon^{-2}\|\nabla^k(\varphi,\zeta)\|^2\right),\\
I_k^{33} &\overset{(1.18)}{\le} C\left(\left\|\frac{\mathbf{D}(\mathbf{u})}{\rho e_\theta(\rho,\theta)}\right\|_{L^3}\|\nabla^k\mathbf{D}(\mathbf{u})\|+\|\mathbf{D}(\mathbf{u})\|_{L^3}\left\|\nabla^k\left(\frac{\mathbf{D}(\mathbf{u})}{\rho e_\theta(\rho,\theta)}\right)\right\|\right)\|\nabla^k\zeta\|_{L^6}\\
&\quad+\left(\left\|\frac{\mathrm{div}\mathbf{u}}{\rho e_\theta(\rho,\theta)}\right\|_{L^3}\|\nabla^k\mathrm{div}\mathbf{u}\|+\|\mathrm{div}\mathbf{u}\|_{L^3}\left\|\nabla^k\left(\frac{\mathrm{div}\mathbf{u}}{\rho e_\theta(\rho,\theta)}\right)\right\|\right)\|\nabla^k\zeta\|_{L^6}\\
&\overset{(3.6),(1.17)}{\le} C\delta_1\left(\|\nabla^{k+1}\mathbf{u}\|^2+\epsilon^{-2}\|\nabla^{k+1}\zeta\|^2+\epsilon^{-2}\|\nabla^k(\varphi,\zeta)\|^2\right),
\end{aligned}
$$

$$
\begin{aligned}
I_k^{34} &\overset{(1.18)}{\leq} C\epsilon^{-2}\left\|h_9(\varphi,\zeta)\right\|_{L^\infty}\left\|\nabla^k(4\tilde{\sigma}\bar{\theta}^3\zeta-\sigma_a\vartheta)\right\|\left\|\nabla^k\zeta\right\| \\
&\quad + C\epsilon^{-2}\left\|(4\tilde{\sigma}\bar{\theta}^3\zeta-\sigma_a\vartheta)\right\|_{L^\infty}\left\|\nabla^k h_9(\varphi,\zeta)\right\|\left\|\nabla^k\zeta\right\| \\
&\overset{(3.6),(1.17)}{\leq} C\delta_1\epsilon^{-2}\left(\left\|\nabla^k\varphi\right\|^2+\left\|\nabla^k\zeta\right\|^2+\left\|\nabla^k(4\tilde{\sigma}\bar{\theta}^3\zeta-\sigma_a\vartheta)\right\|^2\right), \\
I_k^{35} &\overset{(1.18)}{\leq} C\epsilon^{-2}\left(\left\|h_9(\varphi,\zeta)\zeta\right\|_{L^3}\left\|\nabla^k h_5(\zeta)\right\|+\left\|h_5(\zeta)\right\|_{L^3}\left\|\nabla^k[h_9(\varphi,\zeta)\zeta]\right\|\right)\left\|\nabla^k\zeta\right\|_{L^6} \\
&\overset{(3.6),(1.17)}{\leq} C\delta_1\epsilon^{-2}\left(\left\|\nabla^k(\varphi,\zeta)\right\|^2+\left\|\nabla^{k+1}\zeta\right\|^2\right).
\end{aligned}
$$

Therefore, we have

$$
I_k^3 \leq C\delta_1\left(\left\|\nabla^{k+1}\mathbf{u}\right\|^2+\left\|\nabla^k\mathbf{u}\right\|^2+\epsilon^{-2}\left\|\nabla^{k+1}\zeta\right\|^2+\epsilon^{-2}\left\|\nabla^k(\varphi,\zeta,4\tilde{\sigma}\bar{\theta}^3\zeta-\sigma_a\vartheta)\right\|^2\right). \quad (3.21)
$$

By (3.20) and (3.21), we get

$$
\sum_{k=0}^{l}|I_k^3| \leq C\delta_1\left(\left\|\nabla\mathbf{u}\right\|_{H^l}^2+\epsilon^{-2}\left\|\nabla\varphi\right\|_{H^{l-1}}^2+\epsilon^{-2}\left\|\nabla\zeta\right\|_{H^l}^2+\epsilon^{-2}\left\|(4\tilde{\sigma}\bar{\theta}^3\zeta-\sigma_a\vartheta)\right\|_{H^l}^2\right). \quad (3.22)
$$

For $I_0^4$, using Holder inequality, (3.7) and (1.16), we obtain from (3.14) that

$$
I_0^4 \leq C\epsilon^{-2}\left\|h_5(\zeta)\right\|_{L^3}\left\|\zeta\right\|_{L^6}\left\|(4\tilde{\sigma}\bar{\theta}^3\zeta-\sigma_a\vartheta)\right\| \leq C\epsilon^{-2}\delta_1\left(\left\|\nabla\zeta\right\|^2+\left\|(4\tilde{\sigma}\bar{\theta}^3\zeta-\sigma_a\vartheta)\right\|^2\right).
$$

For $I_k^4(k=1,\cdots,l)$, we obtain from (3.14) that

$$
\begin{aligned}
I_k^4 &= \epsilon^{-2}\int\nabla^k[h_5(\zeta)\zeta]\nabla^k(4\tilde{\sigma}\bar{\theta}^3\zeta-\sigma_a\vartheta)d\mathbf{x} \\
&\leq C\epsilon^{-2}\left(\left\|h_5(\zeta)\right\|_{L^3}\left\|\nabla^k\zeta\right\|_{L^6}+\left\|\zeta\right\|_{L^3}\left\|\nabla h_5(\zeta)\right\|_{L^6}\right)\left\|\nabla^k(4\tilde{\sigma}\bar{\theta}^3\zeta-\sigma_a\vartheta)\right\| \\
&\leq C\delta_1\epsilon^{-2}\left(\left\|\nabla^{k+1}\zeta\right\|^2+\left\|\nabla^k(4\tilde{\sigma}\bar{\theta}^3\zeta-\sigma_a\vartheta)\right\|^2\right).
\end{aligned}
$$

Therefore, we have

$$
\sum_{k=0}^{l}|I_k^4| \leq C\delta_1\epsilon^{-2}\left(\left\|\nabla\zeta\right\|_{H^l}^2+\left\|(4\tilde{\sigma}\bar{\theta}^3\zeta-\sigma_a\vartheta)\right\|_{H^l}^2\right). \quad (3.23)
$$

Substituting (3.16), (3.19), (3.22) and (3.23) into (3.10), and using the smallness of $\delta_1>0$, we get (3.8). The proof of Lemma 3.1 is completed.

Next, we have:

**Lemma 3.2** *Under the assumption* (3.6), *it holds that*

$$
\frac{d}{dt}\sum_{k=0}^{l-1}\int\nabla^k\mathbf{u}\cdot\nabla^{k+1}\varphi d\mathbf{x}+\frac{P_\rho(\bar{\rho},\bar{\theta})}{2\bar{\rho}\epsilon^2}\left\|\nabla\varphi\right\|_{H^{l-1}}^2 \leq C\left\|\nabla\mathbf{u}\right\|_{H^l}^2+C\epsilon^{-2}\left\|\nabla\zeta\right\|_{H^{l-1}}^2 \quad (3.24)
$$

*for any* $t\geq 0$ *and* $l=3,\cdots,N$, *where* $C$ *is the positive constant independent of* $\epsilon,\delta_1$ *and* $t$.

**Proof.** Multiplying $(3.9)_2$ by $\nabla^{k+1}\varphi$ and using

$$
\int\nabla^k\mathbf{u}_t\cdot\nabla^{k+1}\varphi d\mathbf{x}=\frac{d}{dt}\int\nabla^k\mathbf{u}\cdot\nabla^{k+1}\varphi d\mathbf{x}-\bar{\rho}\int(\mathrm{div}\nabla^k\mathbf{u})^2d\mathbf{x}+\int\nabla^k g_1\mathrm{div}\nabla^k\mathbf{u}d\mathbf{x}
$$

due to $(3.9)_1$, we have

$$\frac{d}{dt}\int \nabla^k \mathbf{u}\cdot\nabla^{k+1}\varphi d\mathbf{x} + \frac{P_\rho(\bar\rho,\bar\theta)}{\bar\rho\epsilon^2}\left\|\nabla^{k+1}\varphi\right\|^2 - \bar\rho\int(\mathrm{div}\nabla^k\mathbf{u})^2 d\mathbf{x} = I_k^5 + I_k^6, \tag{3.25}$$

where

$$\begin{aligned} I_k^5 = {} & \epsilon^{-2}\int \nabla^k \mathrm{div}(\varphi\mathbf{u})\mathrm{div}\nabla^k\mathbf{u}d\mathbf{x} - \int\nabla^k[(\mathbf{u},\nabla)\mathbf{u}]\cdot\nabla^{k+1}\varphi d\mathbf{x} \\ & + \epsilon^{-2}\int \nabla^k[h_6(\varphi,\zeta)\nabla\varphi + h_7(\varphi,\zeta)\nabla\zeta]\cdot\nabla^{k+1}\varphi d\mathbf{x} \\ & - \int\nabla^k[h_8(\varphi)(\mu\Delta\mathbf{u} + (\mu+\lambda)\nabla\mathrm{div}\mathbf{u})]\cdot\nabla^{k+1}\varphi d\mathbf{x}, \end{aligned} \tag{3.26}$$

and

$$\begin{aligned} I_k^6 = {} & -\frac{\mu}{\bar\rho}\int\nabla^k\Delta\mathbf{u}\nabla^{k+1}\varphi d\mathbf{x} - \frac{\mu+\lambda}{\bar\rho}\int\nabla^{k+1}\mathrm{div}\mathbf{u}\nabla^{k+1}\varphi d\mathbf{x} \\ & - \frac{P_\theta(\bar\rho,\bar\theta)}{\bar\rho\epsilon^2}\int\nabla^{k+1}\zeta\,\nabla^{k+1}\varphi d\mathbf{x}. \end{aligned} \tag{3.27}$$

Adding (3.25) for $k = 0, 1, \cdots, l-1$ yields that

$$\frac{d}{dt}\sum_{k=0}^{l-1}\int\nabla^k\mathbf{u}\cdot\nabla^{k+1}\varphi d\mathbf{x} + \frac{P_\rho(\bar\rho,\bar\theta)}{2\bar\rho\epsilon^2}\left\|\nabla\varphi\right\|_{H^{l-1}}^2 \le \bar\rho\left\|\nabla\mathbf{u}\right\|_{H^{l-1}}^2 + \sum_{k=0}^{l-1}|I_k^5| + \sum_{k=0}^{l-1}|I_k^6|. \tag{3.28}$$

We estimate the right hand in (3.28). For $I_k^5$, by the same lines as in (3.16) and (3.19), we obtain from (3.26) that

$$\sum_{k=0}^{l-1}|I_k^5| \le C\delta_1\left(\left\|\nabla\mathbf{u}\right\|_{H^l}^2 + \epsilon^{-2}\left\|\nabla\varphi\right\|_{H^{l-1}}^2 + \epsilon^{-2}\left\|\nabla\zeta\right\|_{H^{l-1}}^2\right). \tag{3.29}$$

For $I_k^6$, we obtain from (3.27) that

$$\begin{aligned} \sum_{k=0}^{l-1}|I_k^6| & \le C\sum_{k=0}^{l-1}\left(\left\|\nabla^{k+2}\mathbf{u}\right\|\left\|\nabla^{k+1}\varphi\right\| + \epsilon^{-2}\left\|\nabla^{k+1}\zeta\right\|\left\|\nabla^{k+1}\varphi\right\|\right) \\ & \le \eta\epsilon^{-2}\left\|\nabla\varphi\right\|_{H^{l-1}}^2 + \eta^{-1}\left(\left\|\nabla\mathbf{u}\right\|_{H^l}^2 + \epsilon^{-2}\left\|\nabla\zeta\right\|_{H^{l-1}}^2\right) \end{aligned} \tag{3.30}$$

for any $\eta > 0$.

Substituting (3.29)-(3.30) into (3.25), and using the smallness of $\delta_1 > 0$ and $\eta > 0$, we get (3.24). The proof of Lemma 3.2 is completed.

## 3.2 Proof of Theorem 1.2

We first prove the global existence of the smooth solutions to the reformulated system (1.3), (1.4). To this end, we close a priori estimates in the previous subsection. For $N \ge 3$, let $3 \le l \le N$. Then, summing (3.8) and $\beta\times$(3.24), we obtain

$$\begin{aligned} & \frac{d}{dt}\mathrm{E}^l(t) + \left(\frac{\mu}{\bar\rho} - C\beta\right)\left\|\nabla\mathbf{u}\right\|_{H^l}^2 + \left(\frac{\kappa}{\bar\rho\bar\theta} - C\beta\right)\frac{1}{\epsilon^2}\left\|\nabla\zeta\right\|_{H^l}^2 \\ & + \left(\beta\frac{P_\rho(\bar\rho,\bar\theta)}{2\bar\rho} - C\delta_1\right)\frac{1}{\epsilon^2}\left\|\nabla\varphi\right\|_{H^{l-1}}^2 + \frac{\nu\sigma_a}{4\tilde\sigma\bar\rho\bar\theta^4\epsilon^2}\left\|\nabla\vartheta\right\|_{H^l}^2 \le 0, \end{aligned} \tag{3.31}$$

where $\beta \in (0,1]$ is a small constant and

$$
\begin{aligned}
\mathrm{E}^l(t) := \|\mathbf{u}\|_{H^l}^2 + \beta \sum_{k=0}^{l-1} \int \nabla^k \mathbf{u} \cdot \nabla^{k+1} \varphi d\mathbf{x} + \frac{P_\rho(\bar{\rho},\bar{\theta})}{\bar{\rho}^2 \epsilon^2} \|\varphi\|_{H^l}^2 \\
+ \frac{e_\theta(\bar{\rho},\bar{\theta})}{\bar{\theta}\epsilon^2} \|\zeta\|_{H^l}^2 + \frac{\sigma_a}{4\tilde{\sigma}\epsilon\bar{\rho}\bar{\theta}^4} \|\vartheta\|_{H^l}^2 .
\end{aligned}
\tag{3.32}
$$

By (3.32), we can choose the small $\beta \in (0,1]$ and $\delta_1$ independent $\epsilon \in (0,1]$ such that

$$
\mathrm{E}^l(t) \sim \left\| (\mathbf{u}, \epsilon^{-1}\varphi, \epsilon^{-1}\zeta, \epsilon^{-1/2}\vartheta)(t) \right\|_{H^l}^2
$$

uniformly for all $t \geq 0$, and

$$
\frac{\mu}{\bar{\rho}} - C\beta \geq \frac{\mu}{2\bar{\rho}}, \quad \frac{\kappa}{\bar{\rho}\bar{\theta}} - C\beta \geq \frac{\kappa}{2\bar{\rho}\bar{\theta}}, \quad \text{and} \quad \beta \frac{P_\rho(\bar{\rho},\bar{\theta})}{2\bar{\rho}} - C\delta_1 \geq \beta \frac{P_\rho(\bar{\rho},\bar{\theta})}{4\bar{\rho}}. \tag{3.34}
$$

Integrating (3.31) over $t \in [0,\infty)$, and using (3.33), (3.34) and (1.5), we have

$$
\begin{aligned}
&\sup_{t\in\mathbb{R}_+} \|\mathbf{u}(t)\|_{H^l}^2 + \epsilon^{-2} \sup_{t\in\mathbb{R}_+} \|(\varphi,\zeta)(t)\|_{H^l}^2 + \epsilon^{-1} \sup_{t\in\mathbb{R}_+} \|\vartheta(t)\|_{H^l}^2 \\
&+ \int_0^\infty \|\nabla \mathbf{u}^\epsilon\|_{H^l}^2 \, d\tau + \epsilon^{-2} \int_0^\infty \|\nabla\varphi\|_{H^{l-1}}^2 \, d\tau + \epsilon^{-2} \int_0^\infty \|(\nabla\zeta, \nabla\vartheta)\|_{H^l}^2 \, d\tau \\
&\leq C\left( \|\mathbf{u}(0)\|_{H^l}^2 + \epsilon^{-2} \|(\varphi,\zeta)(0)\|_{H^l}^2 + \epsilon^{-1} \|\vartheta(0)\|_{H^l}^2 \right)
\end{aligned}
\tag{3.35}
$$

for $l = 3.\cdots, N$. Therefore, we get

$$
\left\| (\mathbf{u}, \epsilon^{-1}\varphi, \epsilon^{-1}\zeta, \epsilon^{-1/2}\vartheta)(t) \right\|_{H^3}^2 \leq C \left\| (\mathbf{u}, \epsilon^{-1}\varphi, \epsilon^{-1}\zeta, \epsilon^{-1/2}\vartheta)(0) \right\|_{H^3}^2 . \tag{3.36}
$$

Using (1.14) and (3.36), by a standard continuity argument, we can close the a priori estimate (3.6). This in turn allows us to take $l = N$ in (3.35) and we get (1.15).

Now, the proof of existence and uniqueness of global solution is standard. So, we omit the details for brevity. Also, by the same lines as in Subsection 2.3, we can prove the low Mach number limit of the global solutions.

The proof of Theorem 1.2 is completed.

## Conflict of Interest:

The authors declare no conflict of interest.

## Data Availability Statement:

The authors will permit all the data underlying the findings of their manuscripts to be shared by any researchers or groups who are interested in the article.